\def\STable{s}	
\def\draft{n}

\documentclass[12pt]{amsart}
\usepackage{defs,amssymb,dbnsymb,epic,eepic,graphicx,picins}

\theoremstyle{plain}

\newtheorem{theorem}{Theorem}

\newtheorem{proposition}{Proposition}[section]
\newtheorem{lemma}[proposition]{Lemma}
\newtheorem{corollary}[proposition]{Corollary}

\theoremstyle{definition}
\newtheorem{definition}[proposition]{Definition}

\theoremstyle{remark}

\newtheorem{exercise}[proposition]{Exercise}

\newtheorem{remark}[proposition]{Remark}

\newcommand{\prtag}[1]{\tag*{\llap{$#1$\hskip-\displaywidth}}}

\def\printname#1{
  \if\draft y
    \smash{\makebox[0pt]{\hspace{-0.5in}
      \raisebox{8pt}{\tt\tiny #1}}}
  \fi
}

\newcommand{\mathmode}[1]{$#1$}

\newlength{\standardunitlength}
\catcode`\@=11
\long\def\@makecaption#1#2{%
    \vskip 10pt
    \setbox\@tempboxa\hbox{
      \small\sf{\bfcaptionfont #1. }\ignorespaces #2}%
    \ifdim \wd\@tempboxa >\captionwidth {%
        \rightskip=\@captionmargin\leftskip=\@captionmargin
        \unhbox\@tempboxa\par}%
      \else
        \hbox to\hsize{\hfil\box\@tempboxa\hfil}%
    \fi}
\font\bfcaptionfont=cmssbx10 scaled \magstephalf
\newdimen\@captionmargin\@captionmargin=2\parindent
\newdimen\captionwidth\captionwidth=\hsize
\catcode`\@=12

\newcommand{\sign}{\operatorname{sign}}

\def\lbl#1{\label{#1}\printname{#1}}

\def\qed{{\hfill\text{$\Box$}}}

\newlength{\globalparindent}
\newcommand{\udot}{{\mathaccent\cdot\cup}}

\advance \topmargin by -\headheight
\advance \topmargin by -\headsep
\evensidemargin \oddsidemargin
\newcommand{\Zcw}{{{\check{Z}}^\fourwheel}}
\newcommand{\Zc}{\check{Z}}
\newcommand{\Zw}{{Z^\fourwheel}}

\newcommand{\bbQ}{{\mathbb Q}}
\newcommand{\bbR}{{\mathbb R}}
\def\calA{{\mathcal A}}
\def\calB{{\mathcal B}}

\newcommand{\AS}{\text{AS}}
\newcommand{\Aaall}{[\Aa-I--III]}
\newcommand{\Aa}{\text{\sl\AA}}
\newcommand{\Arhus}{\AA rhus}
\newcommand{\IHX}{\text{IHX}}
\newcommand{\LMO}{\text{LMO}}
\newcommand{\STU}{\text{STU}}
\newcommand{\ZLMO}{\hat{Z}^{\mathrm{LMO}}}
\newcommand{\intFG}{\int^{FG}}

\renewcommand{\qed}{~\hfill$\square$}

\newcommand{\unknot}{{\circlearrowleft}}
\newcommand{\OmOm}{{\left\langle\Omega,\Omega\right\rangle}}

\newcommand\tmat[1]{\begin{pmatrix} #1 & -1 \\ 1 & 0 \end{pmatrix}}

\newcommand{\oast}{{\circledast}}
\def\ShacklingElement{{S_{x;x_{1..\ell}}(a_{1..\ell})}}

\allowdisplaybreaks

\begin{document}
\newdimen\captionwidth\captionwidth=\hsize
\ \vspace{-1cm}

\title{A Rational Surgery Formula for the LMO Invariant}

\author{Dror Bar-Natan}
\address{Institute of Mathematics\\
  The Hebrew University\\
  Giv'at-Ram, Jerusalem 91904\\
  Israel
}
\email{drorbn@math.huji.ac.il}
\urladdr{http://www.ma.huji.ac.il/$\sim$drorbn/}

\author{Ruth Lawrence}
\address{Institute of Mathematics\\
  The Hebrew University\\
  Giv'at-Ram, Jerusalem 91904\\
  Israel
}
\email{ruthel@math.huji.ac.il}
\urladdr{http://www.math.lsa.umich.edu/$\sim$ruthjl/}

\thanks{This paper is available electronically from {\tt
  http://www.ma.huji.ac.il/$\sim$drorbn/} and from
  {\tt http://arXiv.org/abs/math.GT/0007045}.
}
\date{
  This edition: Jul.~7,~2000; \quad
  First edition: May 15, 2000; \quad
  Our label: {\tt RationalSurgery}.
}

\begin{abstract}
We write a formula for the \LMO{} invariant of a rational homology sphere
presented as a rational surgery on a link in $S^3$. Our main tool is a
careful use of the \Arhus{} integral and the (now proven) ``Wheels'' and
``Wheeling'' conjectures of B-N, Garoufalidis, Rozansky and Thurston. As
steps, side benefits and asides we give explicit formulas for the values
of the Kontsevich integral on the Hopf link and on Hopf chains, and for
the \LMO{} invariant of lens spaces and Seifert fibered spaces. We find
that the \LMO{} invariant does not separate lens spaces, is far from
separating general Seifert fibered spaces, but does separate Seifert
fibered spaces which are integral homology spheres.
\end{abstract}

\maketitle

\tableofcontents

\section{Introduction}

In their paper ``Wheels, Wheeling, and the Kontsevich Integral of the
Unknot'', B-N, S.~Garoufalidis, L.~Rozansky and
D.~P.~Thurston~\cite{Bar-NatanGaroufalidisRozanskyThurston:WheelsWheeling}
made two conjectures; the ``Wheels'' conjecture about the value of the
Kontsevich integral of the unknot, and the ``Wheeling'' conjecture about
the relationship between the two natural products on the space of
uni-trivalent diagrams. We quote here a paragraph
from~\cite{Bar-NatanGaroufalidisRozanskyThurston:WheelsWheeling},
explaining a part of their motivation:

\begin{quote} {\em
If [the Wheeling conjecture] is true, one would be able to use it along
with [the Wheels conjecture] and known properties of the Kontsevich
integral (such as its behavior under the operations of change of
framing, connected sum, and taking the parallel of a component as
in~\cite{LeMurakami:Parallel}) to get explicit formulas for the
Kontsevich integral of several [\ldots] knots and links. \ldots
Likewise, using [these conjectures] and the hitherto known or
conjectured values of the Kontsevich integral, one would be able to
compute some values of the \LMO{} 3-manifold
invariant~\cite{LeMurakamiOhtsuki:Universal}, using the ``\Arhus{}
integral'' formula of~\Aaall.
} \end{quote}

The ``Wheels'' and ``Wheeling'' conjectures are now theorems
\cite{Kontsevich:DeformationQuantization, Mochizuki:KirillovDuflo,
HinichVaintrob:CyclicOperads, Bar-NatanLeThurston:InPreparation,
ThurstonD:Thesis}, and this seems a good time to proceed with the plan
outlined above. Thus the purpose of our note is to use Wheels and
Wheeling and the \Arhus{} integral to obtain some explicit formulas for
the values of the Kontsevich integral and the \LMO{} invariant on various
simple knots, links and 3-manifolds, as well as a general formula for
the behavior of the \LMO{} invariant under rational surgery over links
(previously such formulas existed only for integral surgery).

\subsection{The main result: a rational surgery formula for the \LMO{}
invariant}

The \LMO{} invariant $\ZLMO(M)$ of a rational homology 3-sphere $M$, which
is given by surgery on some regular integrally framed link $L$ on $S^3$
(we write, $M=S^3_L$; ``regular'' means that the linking matrix of $L$
is non-singular), is a properly normalized version of the integral
of a certain renormalized version of the Kontsevich integral $Z(L)$
of $L$. Thus following~\Aaall, we write (see alternative form at
Equation~\eqref{eq:WheeledAarhus})
\begin{eqnarray}
  \Zc(L) & := & \nu^{\otimes{}X}Z(L)\in\calA(\circlearrowleft_X),
    \lbl{eq:Zcdef} \\
  \Aa_0(L) & := & \int\sigma\Zc(L)\,dX;
    \qquad\qquad\qquad
    Z_\pm := \Aa_0(\circlearrowleft^{\pm 1}) \lbl{eq:AA0} \\
  \ZLMO(M) & := & Z_+^{-\varsigma_+(L)}Z_-^{-\varsigma_-(L)}\Aa_0(L).
    \lbl{eq:LMODef} 
\end{eqnarray}

\noindent In these equations:

\begin{itemize}

\item $X=(x_i)$ is the set of components of $L$,
  $\calA(\circlearrowleft_X)$ is the space of chord diagrams\footnote{
    Throughout this paper all the spaces of diagrams that we will
    consider are graded and we will automatically complete them with
    respect to the grading, so as to allow infinite sums of terms of
    increasing degree.
  }
  whose skeleton consists of $|X|$ directed circles colored by the
  elements of $X$ modulo the usual $4T$, $\STU$, $\AS$ and $\IHX$
  relations, but not divided by the framing independence relation. For
  general background about chord diagrams and related topological and
  Lie algebraic issues see~\cite{Bar-Natan:Vassiliev}.

\item
  $\nu=Z(\circlearrowleft)\in\calA(\circlearrowleft)\cong\calA(\uparrow)$
  is the Kontsevich integral of the zero-framed unknot
  $\circlearrowleft$, regarded as an element of the space
  $\calA(\circlearrowleft)\cong\calA(\uparrow)$ of chord diagrams whose
  skeleton is a single directed circle (modulo the same relations as
  above), or a single directed line (modulo the same).

\item $\nu^{\otimes{}X}$ is the $|X|$'th tensor power of $\nu$, regarded
  as an element of $\calA(\uparrow)^{\otimes X}$, the $X$'th tensor power
  of $\calA(\uparrow)$.  It acts on the Kontsevich integral $Z(L)$ of $L$
  using the usual ``stick in anywhere'' action $\calA(\uparrow)^{\otimes
  X}\otimes\calA(\circlearrowleft_X)\to \calA(\circlearrowleft_X)$.

\item $\sigma\colon\calA(\circlearrowleft_X)\to\calA(\oast_X)$
  is the formal PBW linear isomorphism between the space
  $\calA(\circlearrowleft_X)$ and the space $\calA(\oast_X)$ (denoted
  $\calB^{\text{links}}(X)$ in~\cite{Aarhus:II}) of $X$-marked
  uni-trivalent diagrams modulo $\AS$, $\IHX$ and $X$-flavored link
  relations (see~\cite[Section~5.2]{Aarhus:II}). The map $\sigma$
  is most easily defined as the inverse of the symmetrization
  map $\chi\colon\calA(\oast_X)\to\calA(\circlearrowleft_X)$. If
  $X$ is a singleton, we often suppress it and
  write $\calA(\oast)$, $\calA(\circlearrowleft)$,
  $\sigma\colon\calA(\circlearrowleft)\to\calA(\oast)$, etc. We note
  that $\sigma$ is a homonymous variant of a better known isomorphism
  $\sigma\colon\calA(\uparrow_X\nolinebreak)\to\calA(\ast_X)$, where
  $\calA(\uparrow_X)$ is the same as $\calA(\circlearrowleft_X)$ but
  with the directed circle skeleton components replaced by directed lines
  and $\calA(\ast_X)$ is the standard space of $X$-marked uni-trivalent
  diagrams modulo $\AS$ and $\IHX$, not reduced by link relations.

\item $\int$ is the key ingredient of ``formal integration''.  It can be
  either the $\LMO$-style ``negative dimensional integral'' $\int^{(m)}$,
  or, in the case when $M$ is a rational homology sphere, the ``formal
  Gaussian integration'' $\intFG$ of~\Aaall. (The equality of these
  two integration theories is in~\cite{Aarhus:III}; a sticky leftover
  from~\Aaall, that $\intFG$ is well defined modulo link relations, is
  our Proposition~\ref{prop:sticky}). We note that $\int$ is valued in
  $\calA(\emptyset)$, the space of trivalent diagrams modulo $\AS$ and
  $\IHX$ (that is, unitrivalent diagrams with no univalent vertices), and
  that $\calA(\emptyset)$ is a commutative algebra under disjoint union.

\item $\circlearrowleft^{\pm 1}$ is the $\pm 1$-framed unknot, and
  $\varsigma_+(L)$ and $\varsigma_-(L)$ are respectively, the numbers of
  positive and negative eigenvalues of the covariance matrix of $Z(L)$,
  which is the linking matrix of $L$ (see~\cite[Definition~2.8 and
  Claim~1.10]{Aarhus:I}).

\end{itemize}

It is rather easy to show (see Section~\ref{subsec:TwoStep}) that if a
rational homology 3-sphere $M$ is given by surgery on some rationally
framed link $L$, then its \LMO{} invariant $\ZLMO(M)$ can be computed
using exactly the same formulas \eqref{eq:Zcdef}--\eqref{eq:LMODef},
only replacing the input $Z(L)$ (which is not defined for rationally
framed links) by some extension thereof, which we will also denote by
$Z(L)$. Our main result in this paper is a precise formula for this
``rationally-framed Kontsevich integral''. We start with some definitions.

\begin{definition}
(\cite{Bar-NatanGaroufalidisRozanskyThurston:WheelsWheeling})
Let $\Omega\in\calA(\ast)$ be given by
\begin{equation} \lbl{eq:OmegaDef}
  \Omega=\exp_\udot \sum_{m=1}^\infty b_{2m}\omega_{2m},
\end{equation}
where the modified Bernoulli numbers $b_{2m}$ are defined
by the power series expansion
\begin{equation} \lbl{eq:bdef}
  \sum_{m=0}^\infty b_{2m}x^{2m} = \frac{1}{2}\log\frac{\sinh x/2}{x/2}
\end{equation}
(so that $b_2=1/48$, $b_4=-1/5760$, etc.) and $\omega_{2m}$ is the
$2m$-wheel, the degree $2m$ uni-trivalent diagram made of a $2m$-gon
with $2m$ legs (so that $\omega_2=\twowheel$, $\omega_4=\fourwheel$, \ldots,
with all vertices oriented counterclockwise).
\end{definition}

\begin{definition} (\Aaall) For any element $A\in\calA(\ast_X)$, define
a map
\[ \hat{A}=\partial_A\colon\calA(\ast_X)\longrightarrow\calA(\ast_X) \]
to act on diagrams $B\in\calA(\ast_X)$ by gluing all legs of $A$ to some
subset of legs of $B$ with matching labels. Likewise, define a pairing
\[
  \langle A,B\rangle_X
  := (\partial_A B)_{X\to 0}
  := \left(\parbox{3.1in}{
    sum of all ways of gluing the $x$-marked legs of $A$ to the
    $x$-marked legs of $B$, for all $x\in X$
  }\right).
\]
We note that $\hat{A}=\partial_A$ and $\langle A,\cdot\rangle$ are well
defined even for arguments $B\in\calA(\oast_X)$, provided $\partial_A$
annihilates all $X$-flavored link relations. This is equivalent to
saying that $A$ is invariant with respect to $x$ for all $x\in X$, where
``invariance'' is defined in Figure~\ref{fig:xinv}.
\end{definition}

\begin{figure}
\[ 
  \printname{xinv}
  \setlength{\unitlength}{0.6\standardunitlength}
  \begin{array}{c}  \hspace{-1.7mm}
    \raisebox{-8pt}{\begingroup\makeatletter\ifx\SetFigFont\undefined
\def\x#1#2#3#4#5#6#7\relax{\def\x{#1#2#3#4#5#6}}%
\expandafter\x\fmtname xxxxxx\relax \def\y{splain}%
\ifx\x\y   
\gdef\SetFigFont#1#2#3{%
  \ifnum #1<17\tiny\else \ifnum #1<20\small\else
  \ifnum #1<24\normalsize\else \ifnum #1<29\large\else
  \ifnum #1<34\Large\else \ifnum #1<41\LARGE\else
     \huge\fi\fi\fi\fi\fi\fi
  \csname #3\endcsname}%
\else
\gdef\SetFigFont#1#2#3{\begingroup
  \count@#1\relax \ifnum 25<\count@\count@25\fi
  \def\x{\endgroup\@setsize\SetFigFont{#2pt}}%
  \expandafter\x
    \csname \romannumeral\the\count@ pt\expandafter\endcsname
    \csname @\romannumeral\the\count@ pt\endcsname
  \csname #3\endcsname}%
\fi
\fi\endgroup
{\renewcommand{\dashlinestretch}{30}
\begin{picture}(12617,3366)(0,-10)
\put(12248.000,1318.000){\arc{6120.662}{2.8431}{3.4401}}
\put(9548.000,1318.000){\arc{6120.662}{5.9847}{6.5817}}
\put(3158,2668){\ellipse{3150}{1350}}
\put(683,1273){\ellipse{1350}{1800}}
\path(7658,1723)(7658,1453)
\path(7883,1723)(7883,1453)
\path(8333,1723)(8333,1453)
\path(7433,1453)(8558,1453)(8558,1183)
	(7433,1183)(7433,1453)
\path(7658,1183)(7658,913)
\path(7883,1183)(7883,913)
\path(8333,1183)(8333,913)
\path(7163,1318)(7433,1318)
\path(9683,1993)(9683,1543)
\path(9908,1993)(9908,1543)
\path(10358,1993)(10358,1543)
\path(9413,1768)(9683,1768)
\path(11483,1993)(11483,1543)
\path(11708,1993)(11708,1543)
\path(12158,1993)(12158,1543)
\path(11213,1768)(11708,1768)
\path(11438,1093)(11438,643)
\path(11663,1093)(11663,643)
\path(12113,1093)(12113,643)
\path(11168,868)(12113,868)
\path(2033,1543)(4283,1543)(4283,1093)
	(2033,1093)(2033,1543)
\path(2483,2083)(2483,1543)
\path(2933,1993)(2933,1543)
\path(3833,2038)(3833,1543)
\path(2483,1093)(2483,643)
\path(2933,1093)(2933,643)
\path(3878,1093)(3878,643)
\path(1358,1318)(2033,1318)
\path(2438,373)(2440,371)(2444,368)
	(2452,363)(2462,355)(2475,346)
	(2491,335)(2507,325)(2524,315)
	(2542,306)(2559,298)(2577,292)
	(2597,287)(2618,283)(2634,281)
	(2652,279)(2672,278)(2693,278)
	(2716,278)(2740,278)(2765,278)
	(2791,279)(2818,280)(2845,281)
	(2872,282)(2899,283)(2925,284)
	(2950,285)(2973,285)(2996,285)
	(3016,285)(3035,285)(3053,284)
	(3068,283)(3094,280)(3116,276)
	(3136,271)(3153,264)(3169,257)
	(3183,250)(3193,244)(3200,240)(3203,238)
\path(3968,373)(3966,371)(3962,368)
	(3954,363)(3944,355)(3931,346)
	(3915,335)(3899,325)(3882,315)
	(3864,306)(3847,298)(3829,292)
	(3809,287)(3788,283)(3772,281)
	(3754,279)(3734,278)(3713,278)
	(3690,278)(3666,278)(3641,278)
	(3615,279)(3588,280)(3561,281)
	(3534,282)(3507,283)(3481,284)
	(3456,285)(3433,285)(3410,285)
	(3390,285)(3371,285)(3353,284)
	(3338,283)(3312,280)(3290,276)
	(3270,271)(3253,264)(3237,257)
	(3223,250)(3213,244)(3206,240)(3203,238)
\put(4733,1228){\makebox(0,0)[lb]{\smash{{\mathmode{=\quad 0;}}}}}
\put(6083,1228){\makebox(0,0)[lb]{\smash{{{\rm where}}}}}
\put(8018,1543){\makebox(0,0)[lb]{\smash{{\mathmode{\scriptstyle \cdots}}}}}
\put(8018,1003){\makebox(0,0)[lb]{\smash{{\mathmode{\scriptstyle \cdots}}}}}
\put(8783,1228){\makebox(0,0)[lb]{\smash{{\mathmode{:=}}}}}
\put(10043,1858){\makebox(0,0)[lb]{\smash{{\mathmode{\scriptstyle \cdots}}}}}
\put(11843,1858){\makebox(0,0)[lb]{\smash{{\mathmode{\scriptstyle \cdots}}}}}
\put(11798,958){\makebox(0,0)[lb]{\smash{{\mathmode{\scriptstyle \cdots}}}}}
\put(10763,1723){\makebox(0,0)[lb]{\smash{{\mathmode{+}}}}}
\put(10763,823){\makebox(0,0)[lb]{\smash{{\mathmode{+}}}}}
\put(9413,823){\makebox(0,0)[lb]{\smash{{\mathmode{+}}}}}
\put(9818,823){\makebox(0,0)[lb]{\smash{{\mathmode{\cdots}}}}}
\put(2528,463){\makebox(0,0)[lb]{\smash{{\mathmode{\scriptstyle x}}}}}
\put(2978,463){\makebox(0,0)[lb]{\smash{{\mathmode{\scriptstyle x}}}}}
\put(2978,463){\makebox(0,0)[lb]{\smash{{\mathmode{\scriptstyle x}}}}}
\put(3878,463){\makebox(0,0)[lb]{\smash{{\mathmode{\scriptstyle x}}}}}
\put(3203,1723){\makebox(0,0)[lb]{\smash{{\mathmode{\cdots}}}}}
\put(3203,823){\makebox(0,0)[lb]{\smash{{\mathmode{\cdots}}}}}
\put(2258,58){\makebox(0,0)[lb]{\smash{{{\rm\scriptsize all $x$-marked legs}}}}}
\put(278,913){\makebox(0,0)[lb]{\smash{{{\rm\scriptsize diagram}}}}}
\put(413,1318){\makebox(0,0)[lb]{\smash{{{\rm\scriptsize other}}}}}
\put(503,1723){\makebox(0,0)[lb]{\smash{{{\rm\scriptsize any}}}}}
\put(3068,2623){\makebox(0,0)[lb]{\smash{{\mathmode{A}}}}}
\end{picture}
} }
      \hspace{-1.9mm}
  \end{array}
 \]
\caption{Invariance with respect to $x$ of a diagram $A$.} \lbl{fig:xinv}
\end{figure}

\begin{remark} Strictly speaking, $\partial_A B$ and $\langle
A,B\rangle$ are ill-defined if both $A$ and $B$ contain struts
(see~\cite[Section~2.2]{Aarhus:II}), so that closed vertex-free loops
can be formed by gluing them together. We will not encounter this problem
in this paper.
\end{remark}

\begin{remark} \lbl{rem:PartialNotation} Formally $\partial_A$ acts as if
it were a differential operator on multi-variable polynomials in symbols
indexed by $X$, according to an operator obtained from $A$ by replacing
each label $x$ by $\partial_x$ (see~\cite[Section~2]{Aarhus:II}). Thus
while dealing with specific diagrams we will sometimes use notation
as follows:
\[ 
  \printname{PartialNotation}
  \setlength{\unitlength}{0.3\standardunitlength}
  \begin{array}{c}  \hspace{-1.7mm}
    \raisebox{-8pt}{\begingroup\makeatletter\ifx\SetFigFont\undefined
\def\x#1#2#3#4#5#6#7\relax{\def\x{#1#2#3#4#5#6}}%
\expandafter\x\fmtname xxxxxx\relax \def\y{splain}%
\ifx\x\y   
\gdef\SetFigFont#1#2#3{%
  \ifnum #1<17\tiny\else \ifnum #1<20\small\else
  \ifnum #1<24\normalsize\else \ifnum #1<29\large\else
  \ifnum #1<34\Large\else \ifnum #1<41\LARGE\else
     \huge\fi\fi\fi\fi\fi\fi
  \csname #3\endcsname}%
\else
\gdef\SetFigFont#1#2#3{\begingroup
  \count@#1\relax \ifnum 25<\count@\count@25\fi
  \def\x{\endgroup\@setsize\SetFigFont{#2pt}}%
  \expandafter\x
    \csname \romannumeral\the\count@ pt\expandafter\endcsname
    \csname @\romannumeral\the\count@ pt\endcsname
  \csname #3\endcsname}%
\fi
\fi\endgroup
{\renewcommand{\dashlinestretch}{30}
\begin{picture}(9740,3063)(0,-10)
\path(1650,2229)(2850,2229)(2850,1029)
	(1650,1029)(1650,2229)
\path(1050,2829)(1650,2229)
\path(2850,2229)(3450,2829)
\path(2850,1029)(3450,429)
\path(1650,1029)(1050,429)
\path(7050,2229)(8250,2229)(8250,1029)
	(7050,1029)(7050,2229)
\path(6450,2829)(7050,2229)
\path(8250,2229)(8850,2829)
\path(8250,1029)(8850,429)
\path(7050,1029)(6450,429)
\path(3750,1704)(4650,1704)
\path(3750,1554)(4650,1554)
\path(4425,1929)(4725,1629)(4425,1329)
\put(3450,2904){\makebox(0,0)[lb]{\smash{{\mathmode{\scriptstyle x}}}}}
\put(8925,2904){\makebox(0,0)[lb]{\smash{{\mathmode{\scriptstyle \partial_x}}}}}
\put(825,129){\makebox(0,0)[lb]{\smash{{\mathmode{\scriptstyle x}}}}}
\put(3450,129){\makebox(0,0)[lb]{\smash{{\mathmode{\scriptstyle x}}}}}
\put(6150,54){\makebox(0,0)[lb]{\smash{{\mathmode{\scriptstyle \partial_x}}}}}
\put(8925,54){\makebox(0,0)[lb]{\smash{{\mathmode{\scriptstyle \partial_x}}}}}
\put(5250,1479){\makebox(0,0)[lb]{\smash{{\mathmode{\hat{D}=}}}}}
\put(0,1479){\makebox(0,0)[lb]{\smash{{\mathmode{D=}}}}}
\put(825,2904){\makebox(0,0)[lb]{\smash{{\mathmode{\scriptstyle x}}}}}
\put(6150,2904){\makebox(0,0)[lb]{\smash{{\mathmode{\scriptstyle \partial_x}}}}}
\end{picture}
} }
      \hspace{-1.9mm}
  \end{array}
 \]
\end{remark}

Next we define the {\em Dedekind symbol} $S(p/q)$ of a reduced rational
number $p/q$. A comprehensive source of information about the Dedekind
symbol is~\cite{KirbyMelvin:DedekindSum}, where (among other things)
the well-definededness and equivalence of the definitions below is
discussed.

\begin{definition} \lbl{def:Dedekind}
The {\em Dedekind symbol} $S(p/q)$ of a reduced rational number $p/q$
is defined by the properties
\[
  S(-x)=-S(x), \quad S(x+1)=S(x) \quad\text{and}\quad
  S\left(\frac{p}{q}\right) + S\left(\frac{q}{p}\right)
  = \frac{p}{q} + \frac{q}{p} + \frac{1}{pq} - 3\sign(pq).
\]
Equivalently, $S(p/q)$ is defined by its relation
$S(p/q)=12\sign(q)s(p,q)$ with the Dedekind sum $s(p,q)$, which is given
by either of two formulas
\[ s(p,q)
  := \sum_{k=1}^{|q|-1}
    \left(\!\!\left(\frac{k}{q}\right)\!\!\right)
    \left(\!\!\left(\frac{kp}{q}\right)\!\!\right)
  = \frac{1}{4|q|}\sum_{k=1}^{|q|-1}
    \cot\left(\frac{k\pi}{q}\right)
    \cot\left(\frac{kp\pi}{q}\right),
\]
where $((x))$ is the sawtooth function $((x)):=x-\lfloor x\rfloor-1/2$.
See Table~\ref{tab:STable} on page~\pageref{tab:STable}.
\end{definition}

\begin{definition} \lbl{def:ZLf}
Let $L$ be a rationally framed link, with framing $f_i=p_i/q_i$ on the
component $x_i$ (measured relative to the $0$ framing), and let $L^0$ be
$L$ with all framings replaced by $0$. Set
\[ Z(L) := Z(L^0) \prod_i
  \chi\hat{\Omega}_{x_i}\left(\Omega_{x_i}^{-1}\Omega_{x_i/q_i}\right)
  \exp\left(
    \frac{f_i}{2}\isolatedchord_{x_i}+\frac{S(f_i)}{48}\theta
  \right),
\]
where in this equation:
\begin{itemize}
\item $\Omega^{-1}$ refers to inversion with respect to the disjoint
  union product of $\calA(\ast)$.
\item $\Omega_x^{-1}$ denotes $\Omega^{-1}$ with all of its univalent
  vertices (``legs'') colored $x$.
\item $\Omega_{x/q}$ denotes $\Omega$ with all of its legs colored by
  $x/q$ (meaning that terms with $k$ legs get multiplied by $1/q^k$).
\item $\theta$ denotes the trivalent diagram $
  \printname{theta}
  \setlength{\unitlength}{0.5\standardunitlength}
  \begin{array}{c}  \hspace{-1.7mm}
    \raisebox{-8pt}{%
\begingroup\makeatletter\ifx\SetFigFont\undefined
\def\x#1#2#3#4#5#6#7\relax{\def\x{#1#2#3#4#5#6}}%
\expandafter\x\fmtname xxxxxx\relax \def\y{splain}%
\ifx\x\y   
\gdef\SetFigFont#1#2#3{%
  \ifnum #1<17\tiny\else \ifnum #1<20\small\else
  \ifnum #1<24\normalsize\else \ifnum #1<29\large\else
  \ifnum #1<34\Large\else \ifnum #1<41\LARGE\else
     \huge\fi\fi\fi\fi\fi\fi
  \csname #3\endcsname}%
\else
\gdef\SetFigFont#1#2#3{\begingroup
  \count@#1\relax \ifnum 25<\count@\count@25\fi
  \def\x{\endgroup\@setsize\SetFigFont{#2pt}}%
  \expandafter\x
    \csname \romannumeral\the\count@ pt\expandafter\endcsname
    \csname @\romannumeral\the\count@ pt\endcsname
  \csname #3\endcsname}%
\fi
\fi\endgroup
{\renewcommand{\dashlinestretch}{30}
\begin{picture}(624,629)(0,-10)
\put(312,307){\ellipse{600}{600}}
\path(12,307)(612,307)
\end{picture}
}
 }
      \hspace{-1.9mm}
  \end{array}
$ in
  $\calA(\emptyset)$ and $\isolatedchord_{x_i}$ denotes the ``isolated
  chord'' diagram in $\calA(\uparrow_{x_i})$.
\end{itemize}
\end{definition}

Notice that if all the $f_i$'s are integers, then
$\Omega_{x_i}^{-1}\Omega_{x_i/q_i}=1$ and $S(f_i)=0$, and thus
$Z(L)=Z(L^0)\prod_i e^{f_i\isolatedchord_{x_i}/2}$, and hence the
new definition of $Z$ extends the standard definition of the Kontsevich
integral of an integrally framed link.

We can finally state our main theorem.

\begin{theorem} \lbl{thm:main} (Proof on page~\pageref{pf:main},
alternative formulation in Theorem~\ref{thm:Alt}).
Let $L$ be a rationally framed link and let $M=S^3_{L}$ be
the rational homology 3-sphere obtained from $S^3$ by surgery along
$L$. Then Equations \eqref{eq:Zcdef}--\eqref{eq:LMODef} continue to hold
for the computation of $\ZLMO(M)$, only using Definition~\ref{def:ZLf}
for the definition of $Z(L)$.
\end{theorem}

\begin{remark} The \LMO{} invariant can be generalized to be an invariant
of 3-manifolds with an embedded link, and the surgery formula for the (thus
generalized) \LMO{} invariant holds even if the base manifold is not
necessarily $S^3$. As our proof of Theorem~\ref{thm:main} is completely
local (see below), the theorem generalizes in the obvious manner to the
case when the base manifold is arbitrary and some passive (non-surgery)
link is also present.
\end{remark}

\subsection{Plan of the proof}

It is well known (e.g.~\cite[Section~9H]{Rolfsen:KnotsLinks})
that rational surgery with parameter $p/q$ over a link component
can be achieved by shackling that component with a framed Hopf
chain and then performing integral surgery, in which the framings
$a_{1,\ldots,\ell}$ are related to $p/q$ via a continued fraction
expansion (cf.~\cite[Equation~0.5]{KirbyMelvin:DedekindSum}):

\begin{equation} \lbl{eq:cfrac}
  
  \printname{Shackling}
  \setlength{\unitlength}{0.7\standardunitlength}
  \begin{array}{c}  \hspace{-1.7mm}
    \raisebox{-8pt}{\input draws/Shackling.tex }
      \hspace{-1.9mm}
  \end{array}
 \quad
  \raisebox{8mm}{\hbox{$\displaystyle
    \frac{p}{q} = \langle a_1,\ldots,a_\ell\rangle :=
    -\cfrac{1}{a_1-
      \cfrac{1}{a_2-\ldots
        \cfrac{1}{a_{\ell-1}-
          \frac{1}{a_\ell}
    }}}
  $}}
\end{equation}

Thus we compute the $\LMO$ invariant of $S^3_{L}$ by introducing an
integrally framed shackled version $L^s$ of $L$, in which all of the
original components of $L$ are shackled as above, and then by computing
$\ZLMO(S^3_{L^s})$. We use a lemma, Lemma~\ref{lem:iterative} below, that
allows us to compute the integral in the definition of $\ZLMO(S^3_{L^s})$
in two steps: first we integrate over the Hopf chains, and then over the
original components of $L^0$. The first step is computational in nature,
rather complicated in the technical sense and takes up the bulk of the
proof. When the output of the first step is appropriately normalized,
it turns out to be $Z(L)$, which then needs to be fed into the second
step, which is nothing but a re-run of the procedure in Equations
\eqref{eq:Zcdef}--\eqref{eq:LMODef}. This proves Theorem~\ref{thm:main}.

\subsection{Plan of the paper} In Section~\ref{sec:basic} we discuss
some preliminaries: the Wheels and Wheeling theorems, formal Gaussian
integration and assorted facts regarding continued fractions and
matrices. In Section~\ref{sec:lemmas} we prove some necessary lemmas,
and in Section~\ref{sec:formula} we carry out the computation mentioned
above, of the Kontsevich integral and of surgery over Hopf chains,
and thus finish the proof of Theorem~\ref{thm:main}.

Section~\ref{sec:comp} contains some further computations. Most
importantly, we compute the \LMO{} invariant of arbitrary lens spaces and
of certain Seifert fibered spaces, and find that the \LMO{} invariant
does not separate lens spaces, is far from separating general Seifert
fibered spaces, but it does separate Seifert fibered spaces which are
integral homology spheres.

\subsection{Acknowledgement} We wish to thank R.~Kirby, S.~Levy,
J.~Lieberum, L.~Rozansky, H.~Rubinstein, D.~P.~Thurston and
B.~Weiss for their remarks, suggestions and ideas. In particular
we wish to thank D.~P.~Thurston for his help with the argument in
Proposition~\ref{prop:sticky}. The first author's research at MSRI was
supported in part by NSF grant DMS-9701755. The second author's research
at the Hebrew University was supported in part by a Guastella Fellowship.

\section{Preliminaries} \lbl{sec:basic}

\subsection{Wheels and Wheeling}

Let us start with the statements of the two fundamental theorems, Wheels
and Wheeling.

\begin{theorem}[Wheels,
  \cite{Bar-NatanGaroufalidisRozanskyThurston:WheelsWheeling,
  Bar-NatanLeThurston:InPreparation, ThurstonD:Thesis}] \lbl{thm:Wheels}
The Kontsevich integral of the unknot, $\nu=Z(\circlearrowleft)$,
is the symmetrization of the Wheels element $\Omega$ of
Equation~\eqref{eq:OmegaDef}: $\nu=\chi\Omega$. \qed
\end{theorem}

\begin{theorem}[Wheeling,
  \cite{Bar-NatanGaroufalidisRozanskyThurston:WheelsWheeling,
  Kontsevich:DeformationQuantization, Mochizuki:KirillovDuflo,
  Bar-NatanLeThurston:InPreparation, ThurstonD:Thesis}] \lbl{thm:Wheeling}
The map $\chi\circ\hat{\Omega}_x\colon
\calA(\ast_x)\longrightarrow\calA(\uparrow_x)$
is an algebra isomorphism, so that
\begin{equation} \lbl{eq:Wheeling}
  \chi\hat{\Omega}_x(A\udot{}B)
  = \chi\hat{\Omega}_x(A)\#\chi\hat{\Omega}_x(B),
  \qquad \text{for }A,B\in\calA(\ast_x).
\end{equation}
\qed
\end{theorem}

We will need a slightly more general form of the Wheeling theorem, also
proven in~\cite{Bar-NatanLeThurston:InPreparation, ThurstonD:Thesis}:

\vskip 9pt\par\noindent{\bfseries Theorem \thetheorem'.\hskip .5em}{\itshape
  Equation~\ref{eq:Wheeling} holds even if $A$ and $B$ are allowed to have
  skeleton components and other univalent vertices beyond those labeled
  $x$, provided $A$ and $B$ are both ``invariant with respect to $x$'',
  meaning that they satisfy the relation in Figure~\ref{fig:xinv}.\qed
}\vskip 9pt

We note that Theorems~\ref{thm:Wheels}--\ref{thm:Wheeling}' have a common
generalization as explained and proven
in~\cite{Bar-NatanLeThurston:InPreparation, ThurstonD:Thesis} (see also
Section~\ref{subsec:HopfLink} of this article):

\begin{theorem} \lbl{thm:OpenHopf}
Let $\OpenHopfUp_{\!x}^y$ denote the open Hopf link
with the open component labeled $x$ and the closed component labeled
$y$. Let $\calA(\uparrow_x\!\circledast_y)$ denote the space of
diagrams that have one directed interval skeleton component marked $x$
and (possibly) a number of univalent vertices marked $y$, modulo the
$\STU$, $\AS$ and $\IHX$ relations and modulo $y$-flavored link relations
as in~\cite[Section~5.2]{Aarhus:II}. Then in
$\calA(\uparrow_x\!\circledast_y)$,
\begin{equation*}
  \sigma_yZ(\OpenHopfUp_{\!x}^y)
  = \Omega_y\udot\exp_\#(\botright_x^{\!\!y})
  = \left(1+
      \frac{1}{48}\,y{\text{\large$\twowheel$}}y
      +\ldots
    \right)
    \left(
      \rightarrow_x + \botright_x^{\!\!y}
      + \frac12 \botright_x^{\!\!y}\#\botright_x^{\!\!y}
      + \ldots
    \right)
  \prtag{\qedsymbol}
\end{equation*}
\end{theorem}

\subsection{Formal Gaussian integration}

\begin{definition} \lbl{def:FGI}
Formal Gaussian integration is defined on formal Gaussian-like expressions
(``perturbed $\Lambda$-Gaussians'', for some invertible ``covariance''
$X\times X$ matrix $\Lambda=(l_{ij})$) in $\calA(\ast_XE)$ ($E$ denotes
some arbitrary extra skeleton components) by
\begin{equation} \lbl{eq:FGI}
  \intFG \!\!\! P\exp\left(\frac12 l_{ij}\strutn{x_i}{x_j}\right)dX
  := \left\langle
    \exp\left(-\frac12 l^{ij}\strutu{x_i}{x_j}\right), P
  \right\rangle_X
  \!\!\!= \left.
    \exp\left(-\frac12 l^{ij}\strutu{\partial_{x_i}}{\partial_{x_j}}\right)P
    \right|_{X\to 0}\!\!\!,
\end{equation}
where $\strutn{x_i}{x_j}$ ($\strutu{\partial_{x_i}}{\partial_{x_j}}$)
denotes the $x_ix_j$-colored ($\partial_{x_i}\partial_{x_j}$-colored)
strut, where $P$ is $X$-substantial (involves no struts both of whose ends
are colored by colors in $X$, see~\cite[Definition~2.7]{Aarhus:II}),
and where $(l^{ij})$ is the inverse matrix of $\Lambda=l_{ij}$. (Here
and below we employ the Einstein summation convention).
\end{definition}

Below we will often need to compute the Kontsevich integral of links, and
formal Gaussian integrals thereof. The Kontsevich integral of a link is
valued in a quotient space $\calA(\oast_X)$ of $\calA(\ast_X)$, and thus it
is useful to note that formal Gaussian integration is well defined for
integrands in $\calA(\oast_X)$:

\begin{proposition} \lbl{prop:sticky}
Formal Gaussian integration is well defined on $\calA(\oast_XE)$.
Specifically, if $Z\in\calA(\oast_XE)$ can be written as a perturbed
$\Lambda$-Gaussian in two ways,
\[ Z
  = P_1\exp\left(\frac12l_{ij}\strutn{x_i}{x_j}\right)
  = P_2\exp\left(\frac12l_{ij}\strutn{x_i}{x_j}\right),
\]
where $P_{1,2}$ are $X$-substantial and in $\calA(\ast_XE)$ but the
equality holds in $\calA(\oast_XE)$, then
\[
  \intFG dX\,P_1\exp\left(\frac12l_{ij}\strutn{x_i}{x_j}\right)
  = \intFG dX\,P_2\exp\left(\frac12l_{ij}\strutn{x_i}{x_j}\right)
\]
in $\calA(E)$.
\end{proposition}

\begin{proof} We argue by induction on $|X|$. If $X$ is empty, so is the
statement of the proposition. Otherwise there are two cases.

{\em The lucky case: } If $l_{11}\neq 0$ we can compute the $x_1$ integral
first (\cite[Proposition~2.13]{Aarhus:II}), and we need to show that
\begin{equation} \lbl{eq:x1first}
  \intFG d(X\backslash\{x_1\})\intFG dx_1\,
  (P_1-P_2)\exp\left(\frac12l_{ij}\strutn{x_i}{x_j}\right)
  = 0,
\end{equation}
where we know that
$(P_1-P_2)\exp\left(\frac12l_{ij}\strutn{x_i}{x_j}\right)$ is link relation
equivalent to $0$ via $X$-flavored link relations. Multiplication by
$\strutn{x_1}{x_1}$ is well defined modulo link relations (exercise!),
and thus
\[
  P := (P_1-P_2)\exp\left(\frac12l_{ij}\strutn{x_i}{x_j}\right)
  \exp\left(-\frac12l_{11}\strutn{x_1}{x_1}\right)
\]
is also link relation equivalent to $0$ via $X$-flavored link relations.
By the definition of formal Gaussian integration the inner integral in
Equation~\eqref{eq:x1first} is given by
\[
  \left\langle
    P, \exp\left(-\frac12l_{11}^{-1}\strutu{x_1}{x_1}\right)
  \right\rangle.
\]
The map $D\mapsto\left\langle D,
\exp\left(-\frac12l_{11}^{-1}\strutu{x_1}{x_1}\right)\right\rangle$
kills all $x_1$-flavored link relations (because the
$x_1$-marked strut $\strutu{x_1}{x_1}$ is $x_1$-invariant),
and maps $(X\backslash\{x_1\})$-flavored link relations to
$(X\backslash\{x_1\})$-flavored link relations.  Therefore the
inner integral in Equation~\eqref{eq:x1first} is link equivalent
to $0$ via $(X\backslash\{x_1\})$-flavored link relations. By
the induction hypothesis we now find that the outer integral in
Equation~\eqref{eq:x1first} vanishes.

{\em The ugly case: } If $l_{11}=0$ we consider 
\[ I(\epsilon) := \intFG dX\,
  (P_1-P_2)\exp\left(\frac12l_{ij}\strutn{x_i}{x_j}\right)
  \exp\epsilon\strutn{x_1}{x_1},
\]
where $\epsilon$ is an arbitrary scalar. Multiplication by
$\strutn{x_1}{x_1}$ is well defined modulo link relations, and so the
integrand here remains link equivalent to $0$. Thus by the lucky case,
$I(\epsilon)$ vanishes for all $\epsilon\neq 0$. On the other hand, the
coefficient of every diagram that appears in $I(\epsilon)$ is a rational
function in $\epsilon$ which is non singular at $\epsilon=0$ because
$\Lambda$ is regular. Thus it must be that $I(0)=0$.
\end{proof}

\begin{remark}

The above proof gives us an opportunity to whine and complain about the
state of our understanding of integration in spaces of diagrams. Two
such integration theories exist. The \LMO{} integration theory as
in~\cite{LeMurakamiOhtsuki:Universal, Le:LMO}, and the formal Gaussian
theory of~\Aaall. Neither one of them is satisfactory:

\begin{itemize}

\item Formal Gaussian integration has a solid conceptual foundation;
it is the diagrammatic analogue of an old time favorite, the theory of
Gaussian integrals and Feynman diagrams. But it is defined only for a
restricted kind of integrands, and thus it can only be used to define an
invariant of a restricted class of 3-manifolds, namely rational homology
spheres. And certain aspects of it, such as its being well-defined modulo
link relations (as above), are somewhat tricky.

\item The \LMO{} integration theory is always defined and there's
no problem showing that it is well-defined modulo link relations,
and thus it is superior to formal Gaussian integration, at least in
some sense. But we (the authors) lack a conceptual understanding of
what it really means. It involves the introduction out of thin air of
some strange relations, and one needs to spend some time verifying
that under these new relations the theory does not collapse to
nothing. In~\cite{LeMurakamiOhtsuki:Universal, Le:LMO}, these relations
and the construction as a whole are not interpreted.  In~\cite{Aarhus:III}
the relations are given a semi-satisfactory interpretation in terms of
``negative dimensions''. But whereas formal Gaussian integration is
clearly the diagrammatic counterpart of Gaussian integration over Lie
algebras, the \LMO{} integration theory is not fully understood as the
diagrammatic counterpart of anything (be it integration or anything
else). The situation clearly needs to be resolved. Is the \LMO{} theory
a diagrammatization of something known? What is it? If not, then it
is a genuinely new piece of mathematics. Genuinely new mathematics is
wonderful, but it is a rare commodity. Is the \LMO{} theory really new?

\end{itemize}

In~\cite{Aarhus:III} it is shown that the two integration theories
agree whenever the weaker one is defined. Thus we could have deduced the
previous proposition from the corresponding one for the \LMO{} theory,
which is easier to prove. But then we would have had to rely on a non
trivial theory which is not yet properly understood.

\end{remark}

\subsection{Surgery, continued fractions and matrices} \lbl{subsec:surgery}

The continued fraction expansion in Equation~\eqref{eq:cfrac} has a
slightly refined version
\begin{equation} \lbl{eq:sl2Z}
  \begin{pmatrix} p \\ q \end{pmatrix} =
  \tmat{0}\tmat{a_1}\tmat{a_2}\cdots\tmat{a_\ell}
  \begin{pmatrix} 1 \\ 0 \end{pmatrix},
\end{equation}
in which the signs of $p$ and $q$ are uniquely determined from
$a_{1,\ldots,\ell}$. For convenience and without loss of generality, we
assume that the signs of $p$ and $q$ are fixed so that~\eqref{eq:sl2Z}
holds. This done, we define the integers $u$ and $v$ from the equality
\[
  \begin{pmatrix} p & u \\ q & v \end{pmatrix}
  = A(a_1,\ldots,a_\ell)
  := \tmat{0}\tmat{a_1}\tmat{a_2}\cdots\tmat{a_\ell}.
\]
(cf.~\cite[Lemma~1.9]{KirbyMelvin:DedekindSum}).

Let $\Lambda=(l_{ij})$ be the (tri-diagonal) linking
matrix of the Hopf chain of Equation~\eqref{eq:cfrac}
(cf.~\cite[Equation~0.6]{KirbyMelvin:DedekindSum}),
\[
  \Lambda = \Lambda(a_1,\ldots,a_\ell) := \begin{pmatrix}
    a_1 & 1 & & & \\
    1 & a_2 & 1 \\
    & 1 & \ddots & & & \\
    & & & \ddots & 1 \\
    & & & 1 & a_\ell
  \end{pmatrix}.
\]

\begin{proposition} \lbl{prop:Lambda} The four corners of the inverse
matrix $\Lambda^{-1}=(l^{ij})$ of $\Lambda$ are given by
\begin{equation} \lbl{eq:LambdaInverse}
  l^{11} = -\frac{p}{q}, \qquad
  l^{1\ell} = l^{\ell 1} = \frac{(-1)^\ell}{q}, \qquad
  l^{\ell\ell} = -\frac{v}{q}.
\end{equation}
\end{proposition}

\begin{proof}
By induction on $\ell$ (and row expansion of the relevant determinants)
one establishes the equality
\[ A(a_1,\ldots,a_\ell) =
  \begin{pmatrix}
    -\det\Lambda(a_2,\ldots,a_\ell) & \det\Lambda(a_2,\ldots,a_{\ell-1}) \\
    \det\Lambda(a_1,\ldots,a_\ell) & -\det\Lambda(a_1,\ldots,a_{\ell-1}) \\
  \end{pmatrix}.
\]
After that, the theorem follows from simple observations regarding the
determinant $\det\Lambda$ and the minors $\Lambda^{(ij)}$ of the matrix
$\Lambda$. Namely, that $\det\Lambda=\det\Lambda(a_1,\ldots,a_\ell)$, that
$\Lambda^{(11)}$ is $\det\Lambda(a_2,\ldots,a_\ell)$, that
$\Lambda^{(1\ell)}$ is triangular with ones on the diagonal, and that
$\Lambda^{(\ell\ell)}$ is $\det\Lambda(a_1,\ldots,a_{\ell-1})$.
\end{proof}

We note that while the matrix $\Lambda$ is not determined by $p/q$, a
certain combination of its trace $\tau$, its signature $\varsigma$ and the
numbers $p/q$ and $v/q$ appearing in the corners of its inverse matrix
does depend only on $p/q$:

\begin{theorem} (cf.~\cite[Equation~0.8]{KirbyMelvin:DedekindSum})
  \lbl{thm:Scombination}
The Dedekind symbol $S(p/q)$ of $p/q$ (see Table~\ref{tab:STable}) is given by
\begin{equation*}
  S(p/q) = 3\varsigma-\tau-l^{11}-l^{\ell\ell} = 3\varsigma-\tau+\frac{p+v}{q}.
  \prtag{\qedsymbol}
\end{equation*}
\end{theorem}

\if\STable s
{\begin{table}
  \[
    \def\negative#1{-#1}
    \begin{array}{ccccccccccccc}
  p & 1 & 2 & 3 & 4 & 5 & 6 & 7 & 8 & 9 & 10 & 11 & 12 \\
  q=1 & 0 & 0 & 0 & 0 & 0 & 0 & 0 & 0 & 0 & 0 & 0 & 0 \\
  q=2 & 0 & \cdot & 0 & \cdot & 0 & \cdot & 0 & \cdot & 0 & \cdot & 0 & \cdot \\
  q=3 & 2 & \negative{2} & \cdot & 2 & \negative{2} & \cdot & 2 & \negative{2} & \cdot & 2 & \negative{2} & \cdot \\
  q=4 & 6 & \cdot & \negative{6} & \cdot & 6 & \cdot & \negative{6} & \cdot & 6 & \cdot & \negative{6} & \cdot \\
  q=5 & 12 & 0 & 0 & \negative{12} & \cdot & 12 & 0 & 0 & \negative{12} & \cdot & 12 & 0 \\
  q=6 & 20 & \cdot & \cdot & \cdot & \negative{20} & \cdot & 20 & \cdot & \cdot & \cdot & \negative{20} & \cdot \\
  q=7 & 30 & 6 & \negative{6} & 6 & \negative{6} & \negative{30} & \cdot & 30 & 6 & \negative{6} & 6 & \negative{6} \\
  q=8 & 42 & \cdot & 6 & \cdot & \negative{6} & \cdot & \negative{42} & \cdot & 42 & \cdot & 6 & \cdot \\
  q=9 & 56 & 16 & \cdot & \negative{16} & 16 & \cdot & \negative{16} & \negative{56} & \cdot & 56 & 16 & \cdot \\
  q=10 & 72 & \cdot & 0 & \cdot & \cdot & \cdot & 0 & \cdot & \negative{72} & \cdot & 72 & \cdot \\
  q=11 & 90 & 30 & 18 & 18 & \negative{30} & 30 & \negative{18} & \negative{18} & \negative{30} & \negative{90} & \cdot & 90 \\
  q=12 & 110 & \cdot & \cdot & \cdot & \negative{2} & \cdot & 2 & \cdot & \cdot & \cdot & \negative{110} & \cdot \\
  q=13 & 132 & 48 & 12 & \negative{12} & 0 & \negative{48} & 48 & 0 & 12 & \negative{12} & \negative{48} & \negative{132} \\
  q=14 & 156 & \cdot & 36 & \cdot & 36 & \cdot & \cdot & \cdot & \negative{36} & \cdot & \negative{36} & \cdot \\
  q=15 & 182 & 70 & \cdot & 38 & \cdot & \cdot & \negative{70} & 70 & \cdot & \cdot & \negative{38} & \cdot \\
  q=16 & 210 & \cdot & 30 & \cdot & \negative{30} & \cdot & \negative{18} & \cdot & 18 & \cdot & 30 & \cdot \\
  q=17 & 240 & 96 & 60 & 0 & 12 & 60 & 12 & \negative{96} & 96 & \negative{12} & \negative{60} & \negative{12} \\
  q=18 & 272 & \cdot & \cdot & \cdot & 16 & \cdot & \negative{16} & \cdot & \cdot & \cdot & 16 & \cdot \\
  q=19 & 306 & 126 & 54 & 66 & 66 & \negative{54} & 18 & \negative{18} & \negative{126} & 126 & 18 & \negative{18} \\
  q=20 & 342 & \cdot & 90 & \cdot & \cdot & \cdot & 90 & \cdot & \negative{42} & \cdot & 42 & \cdot \\
  q=21 & 380 & 160 & \cdot & 20 & \negative{20} & \cdot & \cdot & 16 & \cdot & \negative{160} & 160 & \cdot \\
  q=22 & 420 & \cdot & 84 & \cdot & 36 & \cdot & \negative{84} & \cdot & 36 & \cdot & \cdot & \cdot \\
  q=23 & 462 & 198 & 126 & 102 & 42 & 102 & \negative{6} & 126 & \negative{42} & \negative{6} & \negative{198} & 198 \\
  q=24 & 506 & \cdot & \cdot & \cdot & 106 & \cdot & 38 & \cdot & \cdot & \cdot & \negative{74} & \cdot \\
  q=25 & 552 & 240 & 120 & 48 & \cdot & \negative{48} & 0 & \negative{120} & 48 & \cdot & \negative{48} & \negative{240} \\
\end{array}

  \]
  \caption{
    Some values of $q\cdot S(p/q)$ for relatively prime $p,q$. Recall
    that $S$ is antisymmetric and periodic with period $1$, and hence
    all values of $S(p/q)$ with $q\leq 25$ can be recovered from
    this table.  (Change the first line of {\tt RationalSurgery.tex}
    to {\tt$\backslash$def$\backslash$STable$\{$l$\}$} to print a denser
    table that runs up to $q=65$).
  } \lbl{tab:STable}
\end{table}}
\else
{\begin{table}
  \[ 
    \def\negative#1{\underline{#1}}
    \renewcommand{\baselinestretch}{0.84}
    \tiny \setlength\arraycolsep{1.1pt} \fontsize{3.6}{2}
    \begin{array}{ccccccccccccccccccccccccccccccccc}
  p & 1 & 2 & 3 & 4 & 5 & 6 & 7 & 8 & 9 & 10 & 11 & 12 & 13 & 14 & 15 & 16 & 17 & 18 & 19 & 20 & 21 & 22 & 23 & 24 & 25 & 26 & 27 & 28 & 29 & 30 & 31 & 32 \\
  q=1 & 0 & 0 & 0 & 0 & 0 & 0 & 0 & 0 & 0 & 0 & 0 & 0 & 0 & 0 & 0 & 0 & 0 & 0 & 0 & 0 & 0 & 0 & 0 & 0 & 0 & 0 & 0 & 0 & 0 & 0 & 0 & 0 \\
  q=2 & 0 & \cdot & 0 & \cdot & 0 & \cdot & 0 & \cdot & 0 & \cdot & 0 & \cdot & 0 & \cdot & 0 & \cdot & 0 & \cdot & 0 & \cdot & 0 & \cdot & 0 & \cdot & 0 & \cdot & 0 & \cdot & 0 & \cdot & 0 & \cdot \\
  q=3 & 2 & \negative{2} & \cdot & 2 & \negative{2} & \cdot & 2 & \negative{2} & \cdot & 2 & \negative{2} & \cdot & 2 & \negative{2} & \cdot & 2 & \negative{2} & \cdot & 2 & \negative{2} & \cdot & 2 & \negative{2} & \cdot & 2 & \negative{2} & \cdot & 2 & \negative{2} & \cdot & 2 & \negative{2} \\
  q=4 & 6 & \cdot & \negative{6} & \cdot & 6 & \cdot & \negative{6} & \cdot & 6 & \cdot & \negative{6} & \cdot & 6 & \cdot & \negative{6} & \cdot & 6 & \cdot & \negative{6} & \cdot & 6 & \cdot & \negative{6} & \cdot & 6 & \cdot & \negative{6} & \cdot & 6 & \cdot & \negative{6} & \cdot \\
  q=5 & 12 & 0 & 0 & \negative{12} & \cdot & 12 & 0 & 0 & \negative{12} & \cdot & 12 & 0 & 0 & \negative{12} & \cdot & 12 & 0 & 0 & \negative{12} & \cdot & 12 & 0 & 0 & \negative{12} & \cdot & 12 & 0 & 0 & \negative{12} & \cdot & 12 & 0 \\
  q=6 & 20 & \cdot & \cdot & \cdot & \negative{20} & \cdot & 20 & \cdot & \cdot & \cdot & \negative{20} & \cdot & 20 & \cdot & \cdot & \cdot & \negative{20} & \cdot & 20 & \cdot & \cdot & \cdot & \negative{20} & \cdot & 20 & \cdot & \cdot & \cdot & \negative{20} & \cdot & 20 & \cdot \\
  q=7 & 30 & 6 & \negative{6} & 6 & \negative{6} & \negative{30} & \cdot & 30 & 6 & \negative{6} & 6 & \negative{6} & \negative{30} & \cdot & 30 & 6 & \negative{6} & 6 & \negative{6} & \negative{30} & \cdot & 30 & 6 & \negative{6} & 6 & \negative{6} & \negative{30} & \cdot & 30 & 6 & \negative{6} & 6 \\
  q=8 & 42 & \cdot & 6 & \cdot & \negative{6} & \cdot & \negative{42} & \cdot & 42 & \cdot & 6 & \cdot & \negative{6} & \cdot & \negative{42} & \cdot & 42 & \cdot & 6 & \cdot & \negative{6} & \cdot & \negative{42} & \cdot & 42 & \cdot & 6 & \cdot & \negative{6} & \cdot & \negative{42} & \cdot \\
  q=9 & 56 & 16 & \cdot & \negative{16} & 16 & \cdot & \negative{16} & \negative{56} & \cdot & 56 & 16 & \cdot & \negative{16} & 16 & \cdot & \negative{16} & \negative{56} & \cdot & 56 & 16 & \cdot & \negative{16} & 16 & \cdot & \negative{16} & \negative{56} & \cdot & 56 & 16 & \cdot & \negative{16} & 16 \\
  q=10 & 72 & \cdot & 0 & \cdot & \cdot & \cdot & 0 & \cdot & \negative{72} & \cdot & 72 & \cdot & 0 & \cdot & \cdot & \cdot & 0 & \cdot & \negative{72} & \cdot & 72 & \cdot & 0 & \cdot & \cdot & \cdot & 0 & \cdot & \negative{72} & \cdot & 72 & \cdot \\
  q=11 & 90 & 30 & 18 & 18 & \negative{30} & 30 & \negative{18} & \negative{18} & \negative{30} & \negative{90} & \cdot & 90 & 30 & 18 & 18 & \negative{30} & 30 & \negative{18} & \negative{18} & \negative{30} & \negative{90} & \cdot & 90 & 30 & 18 & 18 & \negative{30} & 30 & \negative{18} & \negative{18} & \negative{30} & \negative{90} \\
  q=12 & 110 & \cdot & \cdot & \cdot & \negative{2} & \cdot & 2 & \cdot & \cdot & \cdot & \negative{110} & \cdot & 110 & \cdot & \cdot & \cdot & \negative{2} & \cdot & 2 & \cdot & \cdot & \cdot & \negative{110} & \cdot & 110 & \cdot & \cdot & \cdot & \negative{2} & \cdot & 2 & \cdot \\
  q=13 & 132 & 48 & 12 & \negative{12} & 0 & \negative{48} & 48 & 0 & 12 & \negative{12} & \negative{48} & \negative{132} & \cdot & 132 & 48 & 12 & \negative{12} & 0 & \negative{48} & 48 & 0 & 12 & \negative{12} & \negative{48} & \negative{132} & \cdot & 132 & 48 & 12 & \negative{12} & 0 & \negative{48} \\
  q=14 & 156 & \cdot & 36 & \cdot & 36 & \cdot & \cdot & \cdot & \negative{36} & \cdot & \negative{36} & \cdot & \negative{156} & \cdot & 156 & \cdot & 36 & \cdot & 36 & \cdot & \cdot & \cdot & \negative{36} & \cdot & \negative{36} & \cdot & \negative{156} & \cdot & 156 & \cdot & 36 & \cdot \\
  q=15 & 182 & 70 & \cdot & 38 & \cdot & \cdot & \negative{70} & 70 & \cdot & \cdot & \negative{38} & \cdot & \negative{70} & \negative{182} & \cdot & 182 & 70 & \cdot & 38 & \cdot & \cdot & \negative{70} & 70 & \cdot & \cdot & \negative{38} & \cdot & \negative{70} & \negative{182} & \cdot & 182 & 70 \\
  q=16 & 210 & \cdot & 30 & \cdot & \negative{30} & \cdot & \negative{18} & \cdot & 18 & \cdot & 30 & \cdot & \negative{30} & \cdot & \negative{210} & \cdot & 210 & \cdot & 30 & \cdot & \negative{30} & \cdot & \negative{18} & \cdot & 18 & \cdot & 30 & \cdot & \negative{30} & \cdot & \negative{210} & \cdot \\
  q=17 & 240 & 96 & 60 & 0 & 12 & 60 & 12 & \negative{96} & 96 & \negative{12} & \negative{60} & \negative{12} & 0 & \negative{60} & \negative{96} & \negative{240} & \cdot & 240 & 96 & 60 & 0 & 12 & 60 & 12 & \negative{96} & 96 & \negative{12} & \negative{60} & \negative{12} & 0 & \negative{60} & \negative{96} \\
  q=18 & 272 & \cdot & \cdot & \cdot & 16 & \cdot & \negative{16} & \cdot & \cdot & \cdot & 16 & \cdot & \negative{16} & \cdot & \cdot & \cdot & \negative{272} & \cdot & 272 & \cdot & \cdot & \cdot & 16 & \cdot & \negative{16} & \cdot & \cdot & \cdot & 16 & \cdot & \negative{16} & \cdot \\
  q=19 & 306 & 126 & 54 & 66 & 66 & \negative{54} & 18 & \negative{18} & \negative{126} & 126 & 18 & \negative{18} & 54 & \negative{66} & \negative{66} & \negative{54} & \negative{126} & \negative{306} & \cdot & 306 & 126 & 54 & 66 & 66 & \negative{54} & 18 & \negative{18} & \negative{126} & 126 & 18 & \negative{18} & 54 \\
  q=20 & 342 & \cdot & 90 & \cdot & \cdot & \cdot & 90 & \cdot & \negative{42} & \cdot & 42 & \cdot & \negative{90} & \cdot & \cdot & \cdot & \negative{90} & \cdot & \negative{342} & \cdot & 342 & \cdot & 90 & \cdot & \cdot & \cdot & 90 & \cdot & \negative{42} & \cdot & 42 & \cdot \\
  q=21 & 380 & 160 & \cdot & 20 & \negative{20} & \cdot & \cdot & 16 & \cdot & \negative{160} & 160 & \cdot & \negative{16} & \cdot & \cdot & 20 & \negative{20} & \cdot & \negative{160} & \negative{380} & \cdot & 380 & 160 & \cdot & 20 & \negative{20} & \cdot & \cdot & 16 & \cdot & \negative{160} & 160 \\
  q=22 & 420 & \cdot & 84 & \cdot & 36 & \cdot & \negative{84} & \cdot & 36 & \cdot & \cdot & \cdot & \negative{36} & \cdot & 84 & \cdot & \negative{36} & \cdot & \negative{84} & \cdot & \negative{420} & \cdot & 420 & \cdot & 84 & \cdot & 36 & \cdot & \negative{84} & \cdot & 36 & \cdot \\
  q=23 & 462 & 198 & 126 & 102 & 42 & 102 & \negative{6} & 126 & \negative{42} & \negative{6} & \negative{198} & 198 & 6 & 42 & \negative{126} & 6 & \negative{102} & \negative{42} & \negative{102} & \negative{126} & \negative{198} & \negative{462} & \cdot & 462 & 198 & 126 & 102 & 42 & 102 & \negative{6} & 126 & \negative{42} \\
  q=24 & 506 & \cdot & \cdot & \cdot & 106 & \cdot & 38 & \cdot & \cdot & \cdot & \negative{74} & \cdot & 74 & \cdot & \cdot & \cdot & \negative{38} & \cdot & \negative{106} & \cdot & \cdot & \cdot & \negative{506} & \cdot & 506 & \cdot & \cdot & \cdot & 106 & \cdot & 38 & \cdot \\
  q=25 & 552 & 240 & 120 & 48 & \cdot & \negative{48} & 0 & \negative{120} & 48 & \cdot & \negative{48} & \negative{240} & 240 & 48 & \cdot & \negative{48} & 120 & 0 & 48 & \cdot & \negative{48} & \negative{120} & \negative{240} & \negative{552} & \cdot & 552 & 240 & 120 & 48 & \cdot & \negative{48} & 0 \\
  q=26 & 600 & \cdot & 168 & \cdot & 0 & \cdot & 48 & \cdot & 168 & \cdot & \negative{48} & \cdot & \cdot & \cdot & 48 & \cdot & \negative{168} & \cdot & \negative{48} & \cdot & 0 & \cdot & \negative{168} & \cdot & \negative{600} & \cdot & 600 & \cdot & 168 & \cdot & 0 & \cdot \\
  q=27 & 650 & 286 & \cdot & 146 & 70 & \cdot & 146 & \negative{2} & \cdot & 2 & 70 & \cdot & \negative{286} & 286 & \cdot & \negative{70} & \negative{2} & \cdot & 2 & \negative{146} & \cdot & \negative{70} & \negative{146} & \cdot & \negative{286} & \negative{650} & \cdot & 650 & 286 & \cdot & 146 & 70 \\
  q=28 & 702 & \cdot & 162 & \cdot & 78 & \cdot & \cdot & \cdot & \negative{162} & \cdot & \negative{78} & \cdot & \negative{114} & \cdot & 114 & \cdot & 78 & \cdot & 162 & \cdot & \cdot & \cdot & \negative{78} & \cdot & \negative{162} & \cdot & \negative{702} & \cdot & 702 & \cdot & 162 & \cdot \\
  q=29 & 756 & 336 & 216 & 84 & 156 & 156 & \negative{84} & 48 & \negative{36} & 216 & 48 & 0 & \negative{36} & \negative{336} & 336 & 36 & 0 & \negative{48} & \negative{216} & 36 & \negative{48} & 84 & \negative{156} & \negative{156} & \negative{84} & \negative{216} & \negative{336} & \negative{756} & \cdot & 756 & 336 & 216 \\
  q=30 & 812 & \cdot & \cdot & \cdot & \cdot & \cdot & 20 & \cdot & \cdot & \cdot & 52 & \cdot & 20 & \cdot & \cdot & \cdot & \negative{20} & \cdot & \negative{52} & \cdot & \cdot & \cdot & \negative{20} & \cdot & \cdot & \cdot & \cdot & \cdot & \negative{812} & \cdot & 812 & \cdot \\
  q=31 & 870 & 390 & 210 & 198 & 30 & \negative{30} & 78 & 198 & 78 & \negative{210} & 90 & \negative{6} & \negative{6} & \negative{90} & \negative{390} & 390 & 90 & 6 & 6 & \negative{90} & 210 & \negative{78} & \negative{198} & \negative{78} & 30 & \negative{30} & \negative{198} & \negative{210} & \negative{390} & \negative{870} & \cdot & 870 \\
  q=32 & 930 & \cdot & 270 & \cdot & 114 & \cdot & 30 & \cdot & \negative{30} & \cdot & 270 & \cdot & 114 & \cdot & \negative{162} & \cdot & 162 & \cdot & \negative{114} & \cdot & \negative{270} & \cdot & 30 & \cdot & \negative{30} & \cdot & \negative{114} & \cdot & \negative{270} & \cdot & \negative{930} & \cdot \\
  q=33 & 992 & 448 & \cdot & 128 & 124 & \cdot & 92 & \negative{128} & \cdot & 20 & \cdot & \cdot & \negative{124} & \negative{92} & \cdot & \negative{448} & 448 & \cdot & 92 & 124 & \cdot & \cdot & \negative{20} & \cdot & 128 & \negative{92} & \cdot & \negative{124} & \negative{128} & \cdot & \negative{448} & \negative{992} \\
  q=34 & 1056 & \cdot & 264 & \cdot & 216 & \cdot & 216 & \cdot & 96 & \cdot & \negative{264} & \cdot & 0 & \cdot & \negative{96} & \cdot & \cdot & \cdot & 96 & \cdot & 0 & \cdot & 264 & \cdot & \negative{96} & \cdot & \negative{216} & \cdot & \negative{216} & \cdot & \negative{264} & \cdot \\
  q=35 & 1122 & 510 & 330 & 258 & \cdot & 222 & \cdot & 30 & 258 & \cdot & \negative{78} & 330 & \negative{30} & \cdot & \cdot & \negative{78} & \negative{510} & 510 & 78 & \cdot & \cdot & 30 & \negative{330} & 78 & \cdot & \negative{258} & \negative{30} & \cdot & \negative{222} & \cdot & \negative{258} & \negative{330} \\
  q=36 & 1190 & \cdot & \cdot & \cdot & 70 & \cdot & \negative{70} & \cdot & \cdot & \cdot & \negative{38} & \cdot & 38 & \cdot & \cdot & \cdot & \negative{218} & \cdot & 218 & \cdot & \cdot & \cdot & \negative{38} & \cdot & 38 & \cdot & \cdot & \cdot & 70 & \cdot & \negative{70} & \cdot \\
  q=37 & 1260 & 576 & 324 & 180 & 168 & 0 & 60 & 96 & \negative{180} & 36 & \negative{36} & \negative{324} & 144 & 96 & 168 & 60 & \negative{144} & \negative{576} & 576 & 144 & \negative{60} & \negative{168} & \negative{96} & \negative{144} & 324 & 36 & \negative{36} & 180 & \negative{96} & \negative{60} & 0 & \negative{168} \\
  q=38 & 1332 & \cdot & 396 & \cdot & 180 & \cdot & 132 & \cdot & \negative{12} & \cdot & 132 & \cdot & 396 & \cdot & \negative{180} & \cdot & \negative{12} & \cdot & \cdot & \cdot & 12 & \cdot & 180 & \cdot & \negative{396} & \cdot & \negative{132} & \cdot & 12 & \cdot & \negative{132} & \cdot \\
  q=39 & 1406 & 646 & \cdot & 326 & 286 & \cdot & 74 & 286 & \cdot & 326 & \negative{74} & \cdot & \cdot & 106 & \cdot & 38 & \negative{38} & \cdot & \negative{646} & 646 & \cdot & 38 & \negative{38} & \cdot & \negative{106} & \cdot & \cdot & 74 & \negative{326} & \cdot & \negative{286} & \negative{74} \\
  q=40 & 1482 & \cdot & 390 & \cdot & \cdot & \cdot & 150 & \cdot & 138 & \cdot & 102 & \cdot & \negative{390} & \cdot & \cdot & \cdot & \negative{150} & \cdot & \negative{282} & \cdot & 282 & \cdot & 150 & \cdot & \cdot & \cdot & 390 & \cdot & \negative{102} & \cdot & \negative{138} & \cdot \\
  q=41 & 1560 & 720 & 468 & 240 & 120 & 300 & 300 & \negative{120} & 0 & \negative{240} & 108 & 36 & \negative{132} & 468 & 108 & \negative{48} & \negative{36} & \negative{48} & \negative{132} & \negative{720} & 720 & 132 & 48 & 36 & 48 & \negative{108} & \negative{468} & 132 & \negative{36} & \negative{108} & 240 & 0 \\
  q=42 & 1640 & \cdot & \cdot & \cdot & 232 & \cdot & \cdot & \cdot & \cdot & \cdot & 160 & \cdot & \negative{16} & \cdot & \cdot & \cdot & 232 & \cdot & \negative{160} & \cdot & \cdot & \cdot & 160 & \cdot & \negative{232} & \cdot & \cdot & \cdot & 16 & \cdot & \negative{160} & \cdot \\
  q=43 & 1722 & 798 & 462 & 402 & 246 & 42 & \negative{42} & 78 & 162 & 66 & 402 & 30 & 66 & \negative{462} & 210 & \negative{78} & \negative{246} & 30 & \negative{162} & \negative{210} & \negative{798} & 798 & 210 & 162 & \negative{30} & 246 & 78 & \negative{210} & 462 & \negative{66} & \negative{30} & \negative{402} \\
  q=44 & 1806 & \cdot & 546 & \cdot & 366 & \cdot & 114 & \cdot & 366 & \cdot & \cdot & \cdot & 30 & \cdot & 546 & \cdot & 30 & \cdot & 114 & \cdot & \negative{354} & \cdot & 354 & \cdot & \negative{114} & \cdot & \negative{30} & \cdot & \negative{546} & \cdot & \negative{30} & \cdot \\
  q=45 & 1892 & 880 & \cdot & 308 & \cdot & \cdot & 200 & 160 & \cdot & \cdot & \negative{308} & \cdot & 200 & \negative{92} & \cdot & 92 & 160 & \cdot & \negative{52} & \cdot & \cdot & \negative{880} & 880 & \cdot & \cdot & 52 & \cdot & \negative{160} & \negative{92} & \cdot & 92 & \negative{200} \\
  q=46 & 1980 & \cdot & 540 & \cdot & 180 & \cdot & 132 & \cdot & \negative{180} & \cdot & \negative{60} & \cdot & \negative{132} & \cdot & \negative{540} & \cdot & 36 & \cdot & 36 & \cdot & \negative{60} & \cdot & \cdot & \cdot & 60 & \cdot & \negative{36} & \cdot & \negative{36} & \cdot & 540 & \cdot \\
  q=47 & 2070 & 966 & 630 & 486 & 306 & 390 & 222 & 390 & 30 & 90 & \negative{6} & 486 & 42 & \negative{90} & \negative{198} & 630 & 6 & \negative{42} & 306 & \negative{222} & 30 & \negative{198} & \negative{966} & 966 & 198 & \negative{30} & 222 & \negative{306} & 42 & \negative{6} & \negative{630} & 198 \\
  q=48 & 2162 & \cdot & \cdot & \cdot & 322 & \cdot & 398 & \cdot & \cdot & \cdot & \negative{2} & \cdot & 2 & \cdot & \cdot & \cdot & 178 & \cdot & \negative{322} & \cdot & \cdot & \cdot & \negative{434} & \cdot & 434 & \cdot & \cdot & \cdot & 322 & \cdot & \negative{178} & \cdot \\
  q=49 & 2256 & 1056 & 624 & 384 & 456 & 96 & \cdot & \negative{96} & 216 & 456 & 216 & \negative{384} & 96 & \cdot & \negative{96} & \negative{624} & 288 & 48 & \negative{48} & 96 & \cdot & \negative{96} & \negative{288} & \negative{1056} & 1056 & 288 & 96 & \cdot & \negative{96} & 48 & \negative{48} & \negative{288} \\
  q=50 & 2352 & \cdot & 720 & \cdot & \cdot & \cdot & 0 & \cdot & 48 & \cdot & \negative{48} & \cdot & 240 & \cdot & \cdot & \cdot & 720 & \cdot & 48 & \cdot & \negative{48} & \cdot & \negative{240} & \cdot & \cdot & \cdot & 240 & \cdot & 48 & \cdot & \negative{48} & \cdot \\
  q=51 & 2450 & 1150 & \cdot & 578 & 250 & \cdot & 182 & 142 & \cdot & \negative{250} & 178 & \cdot & 578 & 178 & \cdot & \negative{70} & \cdot & \cdot & \negative{142} & \negative{110} & \cdot & 182 & \negative{110} & \cdot & \negative{1150} & 1150 & \cdot & 110 & \negative{182} & \cdot & 110 & 142 \\
  q=52 & 2550 & \cdot & 714 & \cdot & 390 & \cdot & 282 & \cdot & 246 & \cdot & 186 & \cdot & \cdot & \cdot & 282 & \cdot & \negative{714} & \cdot & 186 & \cdot & 390 & \cdot & \negative{246} & \cdot & \negative{522} & \cdot & 522 & \cdot & 246 & \cdot & \negative{390} & \cdot \\
  q=53 & 2652 & 1248 & 816 & 468 & 408 & 492 & 204 & 240 & 492 & 132 & 252 & 96 & \negative{468} & 192 & \negative{204} & 132 & \negative{276} & 816 & 192 & 240 & \negative{408} & \negative{96} & 0 & \negative{252} & \negative{276} & \negative{1248} & 1248 & 276 & 252 & 0 & 96 & 408 \\
  q=54 & 2756 & \cdot & \cdot & \cdot & 556 & \cdot & 308 & \cdot & \cdot & \cdot & 556 & \cdot & \negative{124} & \cdot & \cdot & \cdot & \negative{164} & \cdot & 164 & \cdot & \cdot & \cdot & \negative{308} & \cdot & \negative{124} & \cdot & \cdot & \cdot & 124 & \cdot & 308 & \cdot \\
  q=55 & 2862 & 1350 & 810 & 678 & \cdot & 162 & 510 & 510 & \negative{162} & \cdot & \cdot & 90 & 30 & 678 & \cdot & 102 & 30 & \negative{810} & 378 & \cdot & 42 & \cdot & 90 & \negative{102} & \cdot & \negative{378} & \negative{1350} & 1350 & 378 & \cdot & 102 & \negative{90} \\
  q=56 & 2970 & \cdot & 918 & \cdot & 330 & \cdot & \cdot & \cdot & 90 & \cdot & \negative{330} & \cdot & 138 & \cdot & 198 & \cdot & \negative{6} & \cdot & 918 & \cdot & \cdot & \cdot & 6 & \cdot & 90 & \cdot & \negative{618} & \cdot & 618 & \cdot & \negative{90} & \cdot \\
  q=57 & 3080 & 1456 & \cdot & 560 & 484 & \cdot & 56 & \negative{56} & \cdot & 164 & \negative{20} & \cdot & 92 & \negative{560} & \cdot & \negative{16} & \negative{164} & \cdot & \cdot & 268 & \cdot & 92 & 484 & \cdot & \negative{16} & \negative{20} & \cdot & \negative{1456} & 1456 & \cdot & 20 & 16 \\
  q=58 & 3192 & \cdot & 912 & \cdot & 504 & \cdot & 264 & \cdot & 312 & \cdot & 48 & \cdot & 312 & \cdot & 336 & \cdot & 0 & \cdot & \negative{912} & \cdot & \negative{48} & \cdot & \negative{504} & \cdot & 264 & \cdot & \negative{336} & \cdot & \cdot & \cdot & 336 & \cdot \\
  q=59 & 3306 & 1566 & 1026 & 786 & 666 & 606 & 378 & 222 & 114 & 606 & 54 & 666 & \negative{114} & \negative{66} & 786 & \negative{54} & 378 & \negative{18} & \negative{366} & 1026 & 66 & \negative{222} & \negative{18} & 174 & \negative{126} & \negative{126} & \negative{174} & \negative{366} & \negative{1566} & 1566 & 366 & 174 \\
  q=60 & 3422 & \cdot & \cdot & \cdot & \cdot & \cdot & 290 & \cdot & \cdot & \cdot & 322 & \cdot & 110 & \cdot & \cdot & \cdot & \negative{290} & \cdot & \negative{142} & \cdot & \cdot & \cdot & \negative{110} & \cdot & \cdot & \cdot & \cdot & \cdot & \negative{722} & \cdot & 722 & \cdot \\
  q=61 & 3540 & 1680 & 1020 & 660 & 420 & 240 & 408 & 336 & 348 & \negative{240} & 0 & \negative{420} & 60 & \negative{60} & \negative{660} & 180 & 96 & 96 & \negative{180} & \negative{1020} & 480 & 108 & 336 & \negative{192} & 108 & \negative{408} & \negative{348} & \negative{192} & \negative{480} & \negative{1680} & 1680 & 480 \\
  q=62 & 3660 & \cdot & 1140 & \cdot & 588 & \cdot & 636 & \cdot & 636 & \cdot & 276 & \cdot & 180 & \cdot & \negative{204} & \cdot & 276 & \cdot & \negative{180} & \cdot & 1140 & \cdot & \negative{12} & \cdot & 588 & \cdot & \negative{12} & \cdot & \negative{204} & \cdot & \cdot & \cdot \\
  q=63 & 3782 & 1798 & \cdot & 902 & 610 & \cdot & \cdot & 646 & \cdot & 218 & 286 & \cdot & 362 & \cdot & \cdot & 902 & 106 & \cdot & 218 & \negative{254} & \cdot & 254 & 286 & \cdot & \negative{610} & 106 & \cdot & \cdot & \negative{362} & \cdot & \negative{1798} & 1798 \\
  q=64 & 3906 & \cdot & 1134 & \cdot & 786 & \cdot & 126 & \cdot & \negative{126} & \cdot & 366 & \cdot & 786 & \cdot & \negative{66} & \cdot & 66 & \cdot & \negative{18} & \cdot & \negative{1134} & \cdot & 126 & \cdot & \negative{126} & \cdot & \negative{18} & \cdot & \negative{366} & \cdot & \negative{834} & \cdot \\
  q=65 & 4032 & 1920 & 1260 & 768 & \cdot & 732 & 360 & 0 & 168 & \cdot & 732 & 180 & \cdot & 288 & \cdot & \negative{768} & 300 & 0 & 108 & \cdot & \negative{468} & 1260 & 300 & 108 & \cdot & \cdot & \negative{180} & 360 & 168 & \cdot & \negative{468} & \negative{1920} \\
\end{array}

  \]
  \caption{
    Some values of $q\cdot S(p/q)$ for relatively prime $p,q$
    (negative numbers are underlined). Recall that $S$ is
    antisymmetric and periodic with period $1$, and hence all
    values of $S(p/q)$ with $q\leq 65$ can be recovered from this
    table. (Change the first line of {\tt RationalSurgery.tex} to
    {\tt$\backslash$def$\backslash$STable$\{$s$\}$} to print a smaller
    table at normal font size).
  } \lbl{tab:STable}
  \end{table}}
\fi

\section{Some Lemmas and preliminary computations} \lbl{sec:lemmas}

\subsection{Simple diagram manipulations}

In the previous section we have used elements of $\calA(\ast_X)$ where $X$
is a set of labels.  One may extend this notation in a multilinear way
so as to allow labels which are elements of a formal vector space whose
basis is indexed by $X$.  Thinking of the elements of $\calA(\ast_X)$
as multinomials in the symbols $X$, whose coefficients are diagrams, one
can also perform evaluations in which one or more of the symbols in $X$
are set to zero, which has the effect of selecting out those diagrams
not mentioning these symbols.

\begin{lemma} \lbl{lem:diff}
For all $A,B\in\calA(\ast_X)$,
$\partial_A(B)=\langle{A_{\bar{X}},B_{X+\bar{X}}}\rangle_{\bar{X}}$.
Here $\bar{X}$ denotes a set in 1--1 correspondence to $X$, $X+{\bar{X}}$
is the set of corresponding sums, $A_{\bar{X}}$ denotes the same diagram
as $A$ with all labels changed from $X$ labels to their corresponding
${\bar{X}}$-label, while the pairing takes place over the set of labels
$\bar{X}$. \qed
\end{lemma}

\begin{lemma} \lbl{lem:DDelta}
Let $E$ be some skeleton specification (a list of colored $\uparrow$'s,
$\circlearrowleft$'s, $\ast$'s and $\circledast$'s) not containing
the colors $t$, $x$ and $y$.  For any $D\in\calA(\ast_t)$ and
$A\in\calA(\ast_t E)$,
\begin{equation*}
  \left[\smash{\hat{D}_t}A\right]_{t\to x+y}
  = \hat{D}_x\left[A\right]_{t\to x+y}
  = \hat{D}_y\left[A\right]_{t\to x+y}.
  \prtag{\qedsymbol}
\end{equation*}
\end{lemma}

\subsection{Properties of $\Omega$}

\begin{proposition}[Pseudo-linearity of $\log\Omega$] \lbl{prop:PseudoLinearity}
For any $x$-invariant diagram $D$,
\begin{equation} \lbl{eq:PseudoLinearity}
  \partial_D\Omega_x
  = \left.\partial_D\Omega_x\right|_{x=0}\Omega_x
  = \langle D,\Omega_x\rangle_x\Omega_x
\end{equation}
(compare with standard calculus: if $D$ is any differential operator and
$f$ is a linear function, then $De^f=(Df)(0)e^f$; we added the prefix
``pseudo'' above because Equation~\ref{eq:PseudoLinearity} does not hold
for every $D$, but only for $x$-invariant $D$'s).
\end{proposition}

\begin{proof} In~\cite{Bar-NatanLeThurston:InPreparation,
ThurstonD:Thesis} it is shown that $\Omega_{x+y}=\Omega_x\udot\Omega_y$ in
$\calA(\circledast_x\circledast_y)$, and thus using Lemma~\ref{lem:diff},
\[ \partial_D(\Omega)
  = \langle{}D_y,\Omega_{x+y}\rangle_y
  = \langle{}D_y,\Omega_x\udot\Omega_y\rangle_y
  = \langle{}D_y,\Omega_y\rangle_y\Omega_x
  = \langle{}D_x,\Omega_x\rangle_x\Omega_x\,.
\]
The $x$-invariance of $D$ is used in asserting the equality
between the contractions $\langle{}D_y,\Omega_{x+y}\rangle_y$
and $\langle{}D_y,\Omega_x\udot\Omega_y\rangle_y$. If
$x$-invariance is not assumed, the contraction map $\langle
D_y,\cdot\rangle_y\colon\calA(\ast_x\ast_y)\to\calA(\ast_x)$
may not descend to a map
$\calA(\circledast_x\circledast_y)\to\calA(\circledast_x)$.
\end{proof}

\begin{corollary} \lbl{cor:ExpDOmega} For any $D\in\calA(\ast)$ we have
$e^{\partial_D}\Omega=e^{\langle D,\Omega\rangle}\Omega$.  In particular,
for any $\alpha$,
\[
  \exp\left(\frac{\alpha}{2}\strutu{\partial_x}{\partial_x}\right)\Omega_x
  = \exp\left(\frac{\alpha\theta}{48}\right)\Omega_x
\]
\end{corollary}

\begin{proof} The first assertion follows immediately
from Proposition~\ref{prop:PseudoLinearity}. The second
assertion follows from the first and from the equality
$\langle\strutu{}{},\Omega\rangle=\theta/24$, which follows from the
fact that the two-legged part of $\Omega$ is $\twowheel/48$.
\end{proof}

\begin{corollary} \lbl{cor:multiplicativity} For any $A,B\in\calA(\ast)$,
\[
  \langle A\udot B,\Omega\rangle
  = \langle A,\Omega\rangle\langle B,\Omega\rangle.
\]
In particular, if $A$ is invertible, then
$\langle{}A^{-1},\Omega\rangle=\langle{}A,\Omega\rangle^{-1}$.
\end{corollary}

\begin{proof}
Clearly, $\partial_{A\udot{}B}=\partial_A\circ\partial_B$. And so,
using Proposition~\ref{prop:PseudoLinearity} twice and the fact that
$\left.\Omega\right|_{x=0}=1$, we get
\begin{align*}
  \langle A\udot B,\Omega\rangle
  &= \left.\partial_{A\udot{}B}\Omega\right|_{x=0}
  = \left.\partial_A(\partial_B(\Omega))\right|_{x=0} \\
  &= \left.\partial_A(\Omega)\right|_{x=0}\cdot\langle B,\Omega\rangle
  = \left.\Omega\right|_{x=0}
    \cdot\langle A,\Omega\rangle\langle B,\Omega\rangle
  = \langle A,\Omega\rangle\langle B,\Omega\rangle.
  \prtag{\qedsymbol}
\end{align*}
\renewcommand\qed{}
\end{proof}

\begin{lemma} \lbl{lem:WiseGuess} The following two formulas hold when
observed\footnote{See~\cite{Mermin:Moon}.}:
\begin{eqnarray}
  \lbl{eq:WiseGuess1}
    \hat{\Omega}\exp\left(\frac{\strutn{}{}}{2}\right)
    & = & \exp\left(\frac{\theta}{48}\right)\Omega\,
      \exp\left(\frac{\strutn{}{}}{2}\right) \\
  \lbl{eq:WiseGuess2}
    \hat{\Omega}_x\exp\left(\frac{\strutn{x+y}{x+y}}{2}\right)
    & = & \exp\left(\frac{\theta}{48}\right)\Omega_{x+y}
      \exp\left(\frac{\strutn{x+y}{x+y}}{2}\right)
\end{eqnarray}
\end{lemma}

\begin{proof} Equation~\eqref{eq:WiseGuess1} follows from the following
computation:
\begin{xalignat*}{2}
  \hat{\Omega}_x\exp\left(\frac{\strutn{x}{x}}{2}\right)
  &= \left\langle
      \exp\left(\frac{\strutn{x+y}{x+y}}{2}\right),\Omega_y
    \right\rangle_y &
    \text{by Lemma~\ref{lem:diff}} \\
  &= \left\langle
      e^{\strutn{x}{x}/2}e^{\strutn{x}{y}}e^{\strutn{y}{y}/2}, \Omega_y
    \right\rangle_y 
    = e^{\strutn{x}{x}/2} \left\langle
      e^{\strutn{x}{y}}e^{\strutn{y}{y}/2}, \Omega_y
    \right\rangle_y & \\
  &= e^{\strutn{x}{x}/2} \left\langle
      e^{\strutn{x}{y}}, e^{\strutu{\partial_y}{\partial_y}/2}\Omega_y
    \right\rangle_y &
    \text{see Remark~\ref{rem:PartialNotation}} \\
  &= \exp\left(\frac{\theta}{48}\right) e^{\strutn{x}{x}/2}
    \left\langle e^{\strutn{x}{y}}, \Omega_y\right\rangle_y &
    \text{by Corollary~\ref{cor:ExpDOmega},} \\*
  &= \exp\left(\frac{\theta}{48}\right) e^{\strutn{x}{x}/2}
    \Omega_x
\end{xalignat*}

The operator $\hat{\Omega}$ commutes with translation by $y$ (the map
$x\mapsto x+y$), because $\hat{\Omega}$ has ``constant coefficients''.
Thus Equation~\eqref{eq:WiseGuess2} is a consequence of
Equation~\eqref{eq:WiseGuess1}. Alternatively, it can be proven
directly along the same lines.
\end{proof}

\subsection{Properties of $\protect\Zw$} \lbl{subsec:Zw}

The Wheeling Theorem allows us to introduce a variant $\Zw$ of the
Kontsevich integral $Z$ which has some nice multiplicative properties:

\begin{definition} \lbl{def:Zw} If
$L$ is an $X$-marked link (or tangle), set
$\Zw:=\left(\prod_i\hat{\Omega}^{-1}_{x_i}\circ\sigma_{x_i}\right)\circ
Z$ to be the ``wheeled Kontsevich integral''. As $\sigma$ and
$\hat{\Omega}^{-1}$ are invertible, $\Zw$ carries just as much information
as the original Kontsevich integral. In a completely parallel manner, set
$\Zcw:=\left(\prod_i\hat{\Omega}^{-1}_{x_i}\circ\sigma_{x_i}\right)\circ
\Zc$.
\end{definition}

A first example of a nice multiplicative property of $\Zw$ is the following
lemma:

\begin{lemma} \lbl{lem:FramingChange} If $L$ is a (possibly rationally)
framed link and $L^0$ is the corresponding $0$-framed link, then
\[ \Zw(L) = \Zw(L^0) \prod_i
  \exp\left((S(f_i)-f_i)\frac{\theta}{48}\right)
  \Omega_{x_i}^{-1}\Omega_{x_i/q_i}
  \exp\frac{f_i}{2}\strutn{x_i}{x_i}.
\]
In particular, if $L$ is integrally framed, then simply
\[ \Zw(L) = \Zw(L^0) \prod_i
  \exp\left(-\frac{f_i\theta}{48}\right)
  \exp\frac{f_i}{2}\strutn{x_i}{x_i}.
\]
\end{lemma}

\begin{proof} Indeed,
\begin{xalignat*}{2}
  \Zw(L)
  &= \left(\prod_i\hat{\Omega}^{-1}_{x_i}\sigma_{x_i}\right) Z(L) & \\
  &= \Zw(L^0)\prod_i
    \Omega_{x_i}^{-1}\Omega_{x_i/q_i}\exp\left(
      \frac{f_i}{2}\hat{\Omega}_{x_i}^{-1}\sigma_{x_i}\isolatedchord_{x_i}
      + \frac{S(f_i)}{48}\theta
    \right)
    & \parbox{108pt}{\raggedleft
      by Definition~\ref{def:ZLf} and Wheeling
    } \\*
  &= \Zw(L^0)\prod_i
    \Omega_{x_i}^{-1}\Omega_{x_i/q_i}\exp\left(
      \frac{f_i}{2}\left(\strutn{x_i}{x_i}-\frac{\theta}{24}\right)
      + \frac{S(f_i)}{48}\theta
    \right)
    & \parbox{108pt}{\raggedleft
      as $\sigma\isolatedchord=\strutn{}{}$ and as
      $\Omega=1+\frac{\twowheel}{48}+\ldots$.
    }
\end{xalignat*}      
\end{proof}

Another nice multiplicative property of $\Zw$ is its behavior under the
operation of taking the connected sum of links:

\begin{lemma} \lbl{lem:ConnectedSum}
Let $L_i$ be an $X_i$-marked link and let $x_i$ be some specific colors
in $X_i$, for $i=1,2$, with $X_i$ disjoint sets of colors.  Let $t$ be
a color not in $X_1\cup X_2$, and let $L$ be the connected sum of $L_1$
and $L_2$ along the components $x_1$ and $x_2$, with the ``merged''
component marked $t$. Then
\begin{eqnarray}
  \lbl{eq:ConnectedSum1} \Zw(L)
    & = & \Zw(\circlearrowleft_t)^{-1}
    \left(\Zw(L_1)\left/x_1\to t\right.\right)
    \left(\Zw(L_2)\left/x_2\to t\right.\right) \\*
  \lbl{eq:ConnectedSum2}
    & = & \OmOm\Omega_t^{-1}
    \left(\Zw(L_1)\left/x_1\to t\right.\right)
    \left(\Zw(L_2)\left/x_2\to t\right.\right).
\end{eqnarray}
\end{lemma}

\begin{proof} Equation~\eqref{eq:ConnectedSum1} follows
immediately from the generalized form of the Wheeling Theorem,
Theorem~\ref{thm:Wheeling}', and the known multiplicative property of the
Kontsevich integral $Z$. Equation~\eqref{eq:ConnectedSum2} follows from
Equation~\eqref{eq:ConnectedSum1} using the following lemma, which is of
independent interest.
\end{proof}

\begin{lemma} \lbl{lem:WheeledOmega}
The wheeled Kontsevich integral of the 0-framed unknot is given by
\begin{equation} \lbl{eq:WheeledUnknot}
  \Zw(\circlearrowleft) = \hat\Omega^{-1}\sigma\nu = \OmOm^{-1}\Omega.
\end{equation}
\end{lemma}

\begin{proof}
\begin{xalignat*}{2}
  \Zw(\circlearrowleft) = \hat\Omega^{-1}\sigma\nu 
  &= \hat\Omega^{-1}\Omega &
    \text{by Wheels} \\
  &= \left\langle\Omega^{-1},\Omega\right\rangle\Omega, &
    \text{by Proposition~\ref{prop:PseudoLinearity},} \\*
  &= \OmOm^{-1}\Omega, & \text{by Corollary~\ref{cor:multiplicativity}.}
  \prtag{\qedsymbol}
\end{xalignat*}
\renewcommand\qed{}
\end{proof}

\begin{corollary} \lbl{cor:Zcw} If $L$ is an $X$-marked link, then
\begin{equation*}
  \Zcw(L) = \OmOm^{-|X|}\left(\prod_{x\in X}\Omega_x\right)\Zw(L).
  \prtag{\qedsymbol}
\end{equation*}
\end{corollary}

Recall that the integration by parts formula
of~\cite[Proposition~2.15]{Aarhus:II} implies that for any
$D,Z\in\calA(\ast_x)$ one has $\int\hat{D}Zdx=\int Zdx$ ($\hat{D}$ is
divergence-free). In particular, $\int\Zcw dX=\int\sigma\Zc dX$ and thus:

\begin{corollary} \lbl{cor:AaFromZw} If $L$ is an $X$-marked link (possibly
rationally framed), then
\begin{equation} \lbl{eq:AaFromZw}
  \Aa_0(L) = \OmOm^{-|X|}\int\left(\prod_{x\in X}\Omega_x\right)\Zw(L)dX.
\end{equation}
\qed
\end{corollary}

\begin{exercise} Determine the behavior of $\Zw$ under
\begin{itemize}
\item doubling a component (see~\cite{LeMurakami:Parallel} and
  Proposition~\ref{prop:HopfLink} below).
\item dropping a component.
\end{itemize}
\end{exercise}

\subsection{One specific integral} At several points below we will need the
value of a certain specific integral, which we hereby evaluate:

\begin{lemma} \lbl{lem:SpecificIntegral}
For any scalars $\alpha$, $\beta$ and $\gamma$ the following equality
holds:
\[
  I(\alpha,\beta;\gamma)
  := \int \Omega_{\alpha x} \Omega_{\beta x}
    \exp\left(\frac{\gamma}{2}\strutn{x}{x}\right)dx
  = \exp\left(-\frac{\alpha^2+\beta^2}{48\gamma}\theta\right)
    \langle\Omega_x,\Omega_{\alpha\beta x/\gamma}\rangle_x
\]
\end{lemma}

\begin{proof}
\begin{alignat*}{2}
  I(\alpha,\beta;\gamma)
  &= \left.
      \exp\left(-\frac{1}{2\gamma}\strutu{\partial_u}{\partial_u}\right)
      \Omega_{\alpha u} \Omega_{\beta u}
    \right|_{u=0} &    
  \text{by~\eqref{eq:FGI}, with a dummy $u$} \\
  &= \left.
      \exp\left(-\frac{1}{2\gamma}\strutu{\partial_u}{\partial_u}\right)
      \Omega_{\alpha(u+v)}\Omega_{\beta(u-v)}
    \right|_{u=v=0}
    & \text{ introducing a spurious $v$} \\
  &= \left.
      \exp\left(-\frac{1}{2\gamma}
        \strutu
          {(\alpha\partial_x+\beta\partial_y)}
          {(\alpha\partial_x+\beta\partial_y)}
      \right)
      \Omega_x\Omega_y
    \right|_{x=y=0}
    & \parbox{144pt}{\raggedleft
      by a change of variables,
      {\scriptsize$\left\{\!\!\begin{array}{rcl}
        \alpha(u+v)	& \!\!\to\!\! &	x \\
        \beta(u-v)	& \!\!\to\!\! &	y \\
        \partial_u	& \!\!\to\!\! &	\alpha\partial_x+\beta\partial_y \\
      \end{array}\right.$}
      (see~\cite[Exercise~2.5]{Aarhus:II})
    } \\
  &= \left.
      e^{-\frac{\alpha\beta\strutu{\partial_x}{\partial_y}}{\gamma}}
      e^{-\frac{\alpha^2\strutu{\partial_x}{\partial_x}}{2\gamma}}
      e^{-\frac{\beta^2\strutu{\partial_y}{\partial_y}}{2\gamma}}
      \Omega_x\Omega_y
    \right|_{x=y=0} & \\
  &= \left\langle
      e^{-\frac{\alpha^2\strutu{x}{x}}{2\gamma}}, \Omega_x
    \right\rangle_x
    \left\langle
      e^{-\frac{\beta^2\strutu{y}{y}}{2\gamma}}, \Omega_y
    \right\rangle_y
    \left.
      e^{-\frac{\alpha\beta\strutu{\partial_x}{\partial_y}}{\gamma}}
      \Omega_x\Omega_y
    \right|_{x=y=0} \hspace{-36pt} & \\*
  &= \exp\left(\frac{-\alpha^2\theta}{48\gamma}\right)
    \exp\left(\frac{-\beta^2\theta}{48\gamma}\right)
    \langle\Omega_x,\Omega_{\alpha\beta x/\gamma}\rangle_x
    & \text{by Corollary~\ref{cor:ExpDOmega}}
  \prtag{\qedsymbol}
\end{alignat*}
\renewcommand\qed{}
\end{proof}

\subsection{The framed unknot and an alternative formula for $\ZLMO$}
\lbl{subsec:p1lens}

We now know enough to write a cleaner formula for $\ZLMO$. We start by
computing the Kontsevich integral of the $p$-framed unknot (for an integer
$p$). Using Lemmas~\ref{lem:FramingChange} and~\ref{lem:WheeledOmega}
we get
\begin{equation} \lbl{eq:FramedUnknot}
  \Zw(\circlearrowleft^p) = \OmOm^{-1} \exp\left(-\frac{p\theta}{48}\right)
  \exp\left(\frac{p}{2}\strutn{}{}\right)\cdot\Omega.
\end{equation}

\noindent We can now compute $\Aa_0(\circlearrowleft^p)$:
\begin{xalignat}{2}
  \nonumber
    \Aa_0(\circlearrowleft^p)
    &= \OmOm^{-1}\int\Omega_x\Zw(\circlearrowleft^p_x)dx
    & \text{by Corollary~\ref{cor:AaFromZw}} \\
  \nonumber     
    &= \OmOm^{-2}\exp\left(-\frac{p\theta}{48}\right)
      \int \Omega_x^2\exp\left(\frac{p}{2}\strutn{x}{x}\right)dx &
    \text{by Equation~\eqref{eq:FramedUnknot},} \\
  \nonumber
    &= \OmOm^{-2}
      \exp\left(\frac{-p\theta}{48}\right)
      \exp\left(\frac{-\theta}{24p}\right)
      \langle\Omega_x,\Omega_{x/p}\rangle_x &
    \parbox{96pt}{\raggedleft
      by Lemma~\ref{lem:SpecificIntegral} for $I(1,1;p)$,
    } \\*
  \lbl{eq:PrenormalizedLens}
    &= \langle\Omega_x,\Omega_x^{-2}\Omega_{x/p}\rangle_x
      \exp\left[-\left(\frac{p}{48}+\frac1{24p}\right)\theta\right]\,.
\end{xalignat}

\noindent In particular, the normalization factors in the definition of the
\Arhus{} integral are
\begin{equation} \lbl{eq:normalization}
  Z_\pm = \Aa_0(\circlearrowleft^{\pm1}) = \OmOm^{-1} \exp\frac{\mp\theta}{16}
\end{equation}

We can now write a rather clean formula for $\ZLMO$ in terms of $\Zw$:

\begin{theorem} \lbl{thm:Alt}
Let $L$ be a framed link, let $M=S^3_{L}$ be the result of surgery on $L$,
and let $\varsigma(L)$ be the signature of its linking matrix.
\begin{itemize}
\item If $L$ is integrally framed, then
\begin{equation} \lbl{eq:WheeledAarhus}
  \ZLMO(M) = \exp\left(\frac{\theta}{16}\varsigma(L)\right)
    \int\Zw(L)\prod_i\Omega_{x_i}dx_i.
\end{equation}
\item If $L$ is rationally framed with framing $f_i=p_i/q_i$ on the
component $x_i$ (relative to the $0$ framing) and if $L^0$ is the zero
framed version of $L$, then the following surgery formula is equivalent
to Theorem~\ref{thm:main} (an in particular, as we shall later see,
it holds):
\begin{equation} \lbl{eq:Alt}
  \ZLMO(M) = \exp\left(\frac{\theta}{48}\left(
      3\varsigma(L)+\sum_i S(f_i)-f_i
    \right)\right)
    \int\Zw(L^0)\prod_i e^{f_i\strutn{x_i}{x_i}/2} \Omega_{x_i/q_i}dx_i.
\end{equation}
\end{itemize}
\end{theorem}

\begin{proof} The first assertion, Equation~\eqref{eq:WheeledAarhus} is
a simple assembly of Equations~\eqref{eq:LMODef},~\eqref{eq:AaFromZw}
and~\eqref{eq:normalization}. The second assertion follows by
applying $\hat{\Omega}^{-1}\sigma$ to the definition of $Z(L)$
(Definition~\ref{def:ZLf}), using Wheeling and substituting the result
into Equation~\eqref{eq:WheeledAarhus}.
\end{proof}

\section{The Rational Surgery Formula} \lbl{sec:formula}

\subsection{Proof modulo computations} \lbl{subsec:TwoStep}

Let $L$ be a rationally framed link, with framing $f_i=p_i/q_i$ on the
component $x_i$ (measured relative to the $0$ framing), and let $L^0$
be $L$ with all framing replaced by $0$. Let $L^s$ be a ``shackled''
version of $L$ --- an integrally framed link surgery upon which is
equivalent to surgery upon $L$, as in Section~\ref{subsec:surgery}
and as in Figure~\ref{fig:Ls}. Let $X=\{x_1,\dots,x_n\}$
be the original set of labels for the components of $L$, let
$X^e=x_1^1,\dots,x_1^{\ell_1},\dots,x_n^1,\dots,x_n^{\ell_n}$ be the
set of `extra' labels in $L^s$, and let $X^s=X\udot X^e$ be the full
set of labels of $L^s$. Let $L^e$ be the part of $L^s$ colored by $X^e$
(the Hopf chains).  Let $M=S^3_{L}=S^3_{L^s}$.

\begin{figure}
\[ 
  \printname{Ls}
  \setlength{\unitlength}{0.4\standardunitlength}
  \begin{array}{c}  \hspace{-1.7mm}
    \raisebox{-8pt}{\input draws/Ls.tex }
      \hspace{-1.9mm}
  \end{array}
 \]
\caption{
  The link $L^s$ is obtained from $L$ by changing
  all the framings $p_i/q_i$ to $0$ and shackling all
  the components $x_1,\dots,x_n$ with Hopf chains marked
  $x_1^1,\dots,x_1^{\ell_1},\dots,x_n^1,\dots,x_n^{\ell_n}$ and framed
  $a_1^1,\dots,a_1^{\ell_1},\dots,a_n^1,\dots,a_n^{\ell_n}$, so that
  $p_i/q_i=\langle a_i^1,\dots,a_i^\ell\rangle$ as in
  Equation~\eqref{eq:cfrac}.
} \lbl{fig:Ls}
\end{figure}

In order to prove Theorem~\ref{thm:main}, we compute  $\ZLMO(M)$ in two
steps:
\begin{enumerate}
\item We first run the procedure of Equations
\eqref{eq:Zcdef}--\eqref{eq:LMODef} only on the `extra' components of
$L^s$. Namely, we set
\begin{align}
  \lbl{eq:Step1}
  \check{Z}_1 &:= \nu^{\otimes{}X^e}Z(L^s),
  & \Aa_1 &:= \int\sigma_{X^e}\check{Z}_1\,dX^e,
  & \ZLMO_1 &:= Z_+^{-\varsigma_{1+}}Z_-^{-\varsigma_{1-}}\Aa_1 \\
\intertext{
  (here $\varsigma_{1\pm}$ are the numbers of positive/negative
  eigenvalues of the covariance matrix of $Z(L^s)$ with respect to only
  the variables in $X^e$).
  \item We then run the procedure of Equations
  \eqref{eq:Zcdef}--\eqref{eq:LMODef} over the remaining components,
  taking as our input the result of the first step. Namely, we set
}
  \nonumber
  \check{Z}_2 &:= \nu^{\otimes{}X}\ZLMO_1,
  & \Aa_2 &:= \int\sigma_X\check{Z}_2\,dX,
  & \ZLMO_2 &:= Z_+^{-\varsigma_{2+}}Z_-^{-\varsigma_{2-}}\Aa_2
\end{align}
(here $\varsigma_{1\pm}$ are the numbers of positive/negative eigenvalues of
the covariance matrix of $\ZLMO_1$).
\end{enumerate}

\begin{proof}[Proof of Theorem~\ref{thm:main}] \lbl{pf:main}
Theorem~\ref{thm:main} follows immediately from
Lemma~\ref{lem:iterative} (right below) and Proposition~\ref{prop:main}
(on page~\pageref{prop:main}), which assert that $\ZLMO_2=\ZLMO(M)$
and that $\ZLMO_1=Z(L)$ respectively.
\end{proof}

\begin{lemma} \lbl{lem:iterative} The above two step procedure really
does compute $\ZLMO(M)$. Namely, $\ZLMO_2=\ZLMO(M)$.
\end{lemma}

\begin{proof} This is Theorem~6.6 of~\cite{LeMurakamiOhtsuki:Universal}. It
also follows from the formalism of~\Aaall{} noting that
\begin{itemize}
\item Formal Gaussian integration behaves correctly under iteration
  (\cite[Proposition~2.13]{Aarhus:II}).
\item The covariance matrix of $\ZLMO_1$ is the linking matrix of
  $L$ and if $\varsigma_\pm$ denotes the numbers of positive/negative
  eigenvalues of the covariance matrix of $Z(L^s)$, then
  $\varsigma_{\pm}=\varsigma_{1\pm}+\varsigma_{2\pm}$.
\end{itemize}
\end{proof}

\subsection{The Hopf link and the open Hopf link} \lbl{subsec:HopfLink}

As a first step towards understanding the Kontsevich integral of shackled
links, we compute $\Zw(H(0,0))$, the wheeled Kontsevich integral of the
0-framed Hopf link:

\begin{proposition} \lbl{prop:HopfLink}
The wheeled Kontsevich integral of the $(x,y)$-marked 0-framed Hopf 
link $H_{x,y}(0,0)$ is given by
\[ \Zw(H_{x,y}(0,0)) = \OmOm^{-1} \exp\strutn{x}{y}. \]
\end{proposition}

\begin{proof} The $(1,1)$-framed Hopf link $H_{x,y}(1,1)$ is the
$(x,y)$-marked double $\Delta^t_{x,y}\circlearrowleft^1_t$ of the 1-framed
unknot $\circlearrowleft^1_t$, and hence
\begin{xalignat*}{2}
  \Zw(H_{x,y}(1,1))
  &= \hat{\Omega}_x^{-1}\hat{\Omega}_y^{-1}\sigma_x\sigma_y
    Z(\Delta^t_{x,y}\circlearrowleft^1_t) & \\
  &= \hat{\Omega}_x^{-1}\hat{\Omega}_y^{-1}
    \left[\sigma_tZ(\circlearrowleft^1_t)\right]_{t\to x+y} &
    \text{by~\cite{LeMurakami:Parallel},} \\
  &= \hat{\Omega}_x^{-1}\left[
      \hat{\Omega}_t^{-1}\sigma_tZ(\circlearrowleft^1_t)
    \right]_{t\to x+y} &
    \text{by Lemma~\ref{lem:DDelta},} \\
  &= \hat{\Omega}_x^{-1}\left[
      \Zw(\circlearrowleft^1_t)
    \right]_{t\to x+y} & \\
  &= \OmOm^{-1} \exp\left(-\frac{\theta}{48}\right)
    \hat{\Omega}_x^{-1} \left[
      \exp\left(\frac{1}{2}\strutn{x+y}{x+y}\right)\Omega_{x+y}
    \right] &
    \text{by~\eqref{eq:FramedUnknot},} \\
  &= \OmOm^{-1} \exp\left(-\frac{\theta}{24}\right) 
    \exp\left(\frac{1}{2}\strutn{x+y}{x+y}\right) &
    \text{by Lemma~\ref{lem:WiseGuess},~\eqref{eq:WiseGuess2}.}
\end{xalignat*}
It remains to undo the $(1,1)$ framing by using
Lemma~\ref{lem:FramingChange} on each component:
\begin{alignat*}{1}
  \Zw(H_{x,y}(0,0))
  &= \left[\exp\left(\frac{\theta}{48}\right)\right]^2
    \exp\left(-\frac{1}{2}\strutn{x}{x}\right)
    \exp\left(-\frac{1}{2}\strutn{y}{y}\right)
    \Zw(H_{x,y}(1,1)) \\*
  &= \OmOm^{-1} \exp\strutn{x}{y}.
  \prtag{\qedsymbol}
\end{alignat*}
\renewcommand\qed{}
\end{proof}

For the sake of completeness, we can now re-prove
Theorem~\ref{thm:OpenHopf}:

\begin{proof}[Proof of Theorem~\ref{thm:OpenHopf}]
To get $\sigma_x\sigma_yZ(\OpenHopfUp_{\!x}^y)$ from
$\Zw(H_{x,y}(0,0))$, we need to apply $\hat{\Omega}_x$ and
$\hat{\Omega}_y$, and to ``open up'' the $x$-component, which amounts
to multiplication (using the product $\#$) by $Z(\unknot_x)^{-1}$
(recall that the Kontsevich integral of an open unknot is
trivial). The latter operation can more easily be performed first,
by $\udot$-multiplying $\Zw(H_{x,y}(0,0))$ by $\OmOm\Omega_x^{-1}$,
as in~\eqref{eq:WheeledUnknot}. Thus
\begin{equation} \lbl{eq:WheeledOpenHopf}
  \Zw(\OpenHopfUp_{\!x}^y)
  = \OmOm\Omega_x^{-1}\Zw(H_{x,y}(0,0))
  = \Omega_x^{-1}\exp\strutn{x}{y}
\end{equation}
and hence
\begin{xalignat*}{2}
  \sigma_x\sigma_yZ(\OpenHopfUp_{\!x}^y)
  &= \hat{\Omega}_x\hat{\Omega}_y \Zw(\OpenHopfUp_{\!x}^y) & \\
  &= \hat{\Omega}_x\hat{\Omega}_y
    \left(\Omega_x^{-1}\exp\strutn{x}{y}\right) &
    \text{by Equation~\eqref{eq:WheeledOpenHopf},} \\
  &= \hat{\Omega}_x\left(\Omega_x^{-1}\hat{\Omega}_y\exp\strutn{x}{y}\right)
    = \hat{\Omega}_x\left(\Omega_x^{-1}\Omega_x\exp\strutn{x}{y}\right) & \\*
  &= \hat{\Omega}_x\left(\exp\strutn{x}{y}\right) 
    = \Omega_y\exp\strutn{x}{y}. &
\end{xalignat*}
It only remains to note that
$
  \sigma_yZ(\OpenHopfUp_{\!x}^y)
  = \chi_x\sigma_x\sigma_yZ(\OpenHopfUp_{\!x}^y)
  = \Omega_y\chi_x\exp\strutn{x}{y}
  = \Omega_y\exp_\#(\botright_x^{\!\!y})
$.
\end{proof}

\subsection{The Hopf chain and the shackling element}

We can use the connect-sum lemma (Lemma~\ref{lem:ConnectedSum}) repeatedly
in order to compute the wheeled Kontsevich integral of Hopf chains
(see Figure~\ref{fig:HopfSum}). The result is
\[
  \Zw(H_{x_1,\ldots,x_\ell}(0,\ldots,0)) =
  \OmOm^{-1}\prod_{i=2}^{\ell-1}\Omega_{x_i}^{-1}
  \prod_{i=1}^{\ell-1}\exp\strutn{\ x_i}{x_{i+1}}.
\]
Using Lemma~\ref{lem:FramingChange} and some repackaging, we get:

\begin{figure}
\def\HopfChain{{H_{x_1,\ldots,x_\ell}(a_1,\ldots,a_\ell)}}
\[ 
  \printname{HopfSum}
  \setlength{\unitlength}{0.7\standardunitlength}
  \begin{array}{c}  \hspace{-1.7mm}
    \raisebox{-8pt}{\input draws/HopfSum.tex }
      \hspace{-1.9mm}
  \end{array}
 \]
\caption{    
  A Hopf chain of length $\ell$ is a connected sum of $\ell-1$ Hopf links.
} \lbl{fig:HopfSum}
\end{figure}

\begin{proposition} \lbl{prop:HopfChain} Let $(l_{ij})$ be the linking
matrix of $H_{x_1,\ldots,x_\ell}(a_1,\ldots,a_\ell)$. Then
\begin{equation*}
  \Zw(H_{x_1,\ldots,x_\ell}(a_1,\ldots,a_\ell))
  = \OmOm^{-1}\exp\left(-\frac{\theta\sum a_i}{48}\right)
    \left(\prod_{i=2}^{\ell-1}\Omega_{x_i}^{-1}\right)
    \exp \frac12 l_{ij}\strutn{x_i}{x_j}.
  \prtag{\qedsymbol}
\end{equation*}
\end{proposition}

We can now use Proposition~\ref{prop:HopfChain},
Equation~\eqref{eq:WheeledOpenHopf} and the connect-sum lemma
(Lemma~\ref{lem:ConnectedSum}) to compute the wheeled Kontsevich integral
of the shackling element (see Figure~\ref{fig:ShacklingElement}):

\begin{figure}
\[ 
  \printname{ShacklingElement}
  \setlength{\unitlength}{0.7\standardunitlength}
  \begin{array}{c}  \hspace{-1.7mm}
    \raisebox{-8pt}{\input draws/ShacklingElement.tex }
      \hspace{-1.9mm}
  \end{array}
 \]
\caption{
  The shackling element is a connected sum of an open Hopf link and a Hopf
  chain.
} \lbl{fig:ShacklingElement}
\end{figure}

\begin{proposition} \lbl{prop:ZwS}
Let $(l_{ij})$ be the $\ell\times\ell$ linking matrix of the closed
components of the shackling element $\ShacklingElement$. Then
\begin{equation*}
  \Zw(\ShacklingElement)
  = \Omega_x^{-1} \left(\prod_{i=1}^{\ell-1}\Omega_{x_i}^{-1}\right)
    \exp\left(
      \strutn{\ x}{x_1}
      +\frac12 l_{ij}\strutn{x_i}{x_j}
      - \frac{\theta\sum a_i}{48}
    \right).
  \prtag{\qedsymbol}
\end{equation*}
\end{proposition}

\subsection{The Kontsevich integral of rationally framed links}

\begin{proposition} \lbl{prop:main} The output $\ZLMO_1$ of the first
step of the computation procedure of Section~\ref{subsec:TwoStep} is
equal to the Kontsevich integral for rationally framed links $Z(L)$
of Definition~\ref{def:ZLf}: $\ZLMO_1=Z(L)$.
\end{proposition}

\begin{proof} Both sides of the required equality are clearly
made of local contributions, one per each shackling element
or framing fraction $p_i/q_i$. Thus it is enough to prove the
proposition at the locale where all the actors act. Ergo we may
as well assume that the link $L$ in question is a straight line
$\uparrow^f_x$ marked $f=p/q$, and then, after choosing $\ell$ and
$a_1,\ldots,a_\ell$ as in Equation~\eqref{eq:cfrac} (more precisely,
as in Equation~\eqref{eq:sl2Z}), the shackled $L$ becomes the shackling
element $\ShacklingElement$ (so $X^e$, the set of ``extra'' labels, is
$\{x_1,\dots,x_\ell\}$). We just need to compute $\ZLMO_1$ as in
Section~\ref{subsec:TwoStep}, starting from $L^s=\ShacklingElement$.  Set
$C=\OmOm^{-\ell}\exp\left(-\frac{\theta\sum a_i}{48}\right)\Omega_x^{-1}$
and start crunching:

\begin{alignat}{2}
  \hat{\Omega}^{-1}_x\sigma_x\Aa_1
  &= \OmOm^{-\ell}\int\Zw(\ShacklingElement)
    \prod_i\Omega_{x_i}dx_i
    & \parbox{98pt}{\raggedleft
        by Corollary~\ref{cor:AaFromZw}, applied to $x_{1,\ldots,\ell}$
      } \nonumber \\
  &= C \int\Omega_{x_\ell}\exp\left(
      \strutn{\ x}{x_1}+ \frac12 l_{ij}\strutn{x_i}{x_j}
    \right)\prod_i dx_i
    & \text{by Proposition~\ref{prop:ZwS}} \nonumber \\
  &= C\left\langle
      \exp\left(-\frac12 l^{ij}\strutu{x_i}{x_j}\right),
      \Omega_{x_\ell}e^{\strutn{\ x}{x_1}}
    \right\rangle_{X^e}
    & \text{by Definition~\ref{def:FGI}} \nonumber \\
  &= C\left\langle
      e^{-\frac12l^{11}\strutu{x_1}{x_1}}
      \cdot e^{-l^{1\ell}\strutu{x_1}{x_\ell}}
      \cdot e^{-\frac12l^{\ell\ell}\strutu{x_\ell}{x_\ell}},
      \Omega_{x_\ell}e^{\strutn{\ x}{x_1}}             
    \right\rangle_{x_1,x_\ell}
    & \parbox{98pt}{\raggedleft ignoring inactive variables} \nonumber \\
  &= C\left\langle
      e^{\frac{p}{2q}\strutu{x_1}{x_1}}
      \cdot e^{\frac{-(-1)^{\ell}}{q}\strutu{x_1}{x_\ell}}
      \cdot e^{\frac{v}{2q}\strutu{x_\ell}{x_\ell}},
      \Omega_{x_\ell}e^{\strutn{\ x}{x_1}}
    \right\rangle_{x_1,x_\ell}
    & \text{by Equation~\eqref{eq:LambdaInverse}} \nonumber \\
  &= C\left\langle
      e^{\frac{p}{2q}\strutu{x_1}{x_1}}
      \cdot e^{\frac{-(-1)^{\ell}}{q}\strutu{x_1}{x_\ell}},
      e^{\frac{v}{2q}\strutu{\partial_{x_\ell}}{\partial_{x_\ell}}}
      \Omega_{x_\ell}e^{\strutn{\ x}{x_1}}       
    \right\rangle_{x_1,x_\ell} & \nonumber \\
  &= C \exp\left(\frac{v\theta}{48q}\right) \left\langle
      e^{\frac{p}{2q}\strutu{x_1}{x_1}}
      \cdot e^{\frac{-(-1)^{\ell}}{q}\strutu{x_1}{x_\ell}},
      \Omega_{x_\ell}e^{\strutn{\ x}{x_1}}
    \right\rangle_{x_1,x_\ell}
    & \text{by Corollary~\ref{cor:ExpDOmega}} \nonumber \\
  &= C \exp\left(\frac{v\theta}{48q}\right)
    \Omega_{x/q}\exp\left(\frac{p}{2q}\strutn{x}{x}\right) & \nonumber \\
  &= \OmOm^{-\ell}\exp\left[
      \left(\frac{v}{q}-\sum a_i\right)\frac{\theta}{48}
    \right]
    \Omega_x^{-1}\Omega_{x/q}
    \exp\left(\frac{p}{2q}\strutn{x}{x}\right) & \nonumber \\
  \hat{\Omega}^{-1} \sigma \ZLMO_1
  &= Z_+^{-\varsigma_+} Z_-^{-\varsigma_-} \Aa_1 & \nonumber \\
  &= \exp\left[
      \left(\frac{v}{q}+3\varsigma-\tau\right)\frac{\theta}{48}
    \right]
    \Omega_x^{-1}\Omega_{x/q}
    \exp\left(\frac{p}{2q}\strutn{x}{x}\right)
    & \text{by Equation~\eqref{eq:normalization}} \nonumber \\
  &= \exp\left[
      \left(S(p/q)-p/q\right)\frac{\theta}{48}    
    \right]
    \Omega_x^{-1}\Omega_{x/q}
    \exp\left(\frac{p}{2q}\strutn{x}{x}\right)
    & \text{by Theorem~\ref{thm:Scombination}} \lbl{eq:WheeledFormula} \\
  \ZLMO_1
  &= \exp\left[
      \left(S(p/q)-p/q\right)\frac{\theta}{48}                  
    \right]
    \chi\hat{\Omega}\left(\Omega_x^{-1}\Omega_{x/q}\right)
    \exp\left(\frac{p}{2q}\chi\hat{\Omega}(\strutn{x}{x})\right)
    & \text{by Wheeling} \nonumber \\*
  &= \chi\hat{\Omega}\left(\Omega_x^{-1}\Omega_{x/q}\right)
    \exp\left(\frac{f}{2}\isolatedchord_{x}+\frac{S(f)}{48}\theta\right)
    \nonumber
\end{alignat}
as required.
\end{proof}

\begin{remark} Strictly speaking, our arguments in this section assumed
that $\ell\geq 2$. This is not a serious limitation. One may either read
the arguments again and figure what the appropriate defaults are for
$\ell<2$, or simply use a longer continued fraction expansion for $p/q$.
Notice that we made no restrictive assumptions on the integers $a_i$,
and so every fraction $p/q$ has arbitrarily long continued fraction
expansions.
\end{remark}

\section{Some Computations} \lbl{sec:comp}

\subsection{General $(p,q)$ lens spaces}

As a first application of the rational surgery formula we compute
the \LMO{} invariant of general $(p,q)$ lens spaces.

\begin{proposition} \lbl{prop:ArbitraryLens} The \LMO{} invariant of
arbitrary $(p,q)$ lens spaces is given by
\[
  \ZLMO(L_{p,q})
  = \langle\Omega_x,\Omega_x^{-1}\Omega_{x/p}\rangle_x
    \exp\frac{-S(q/p)}{48}\theta
\]
\end{proposition}

\begin{proof} Recall that the $(p,q)$ lens space $L_{p,q}$
is obtained from $S^3$ by surgery over the $p/q$-framed unknot
$\unknot^{p/q}$. Thus we can use Equation~\eqref{eq:Alt} to compute
$\ZLMO$, remembering also that by Equation~\eqref{eq:WheeledUnknot},
$\Zw(\unknot)=\OmOm^{-1}\Omega$:
\begin{alignat*}{2}
  \ZLMO&(L_{p,q})
  = \OmOm^{-1}\exp\left(\frac{\theta}{48}(3\sign(p/q)+S(p/q)-p/q)\right)
    \int e^{p\strutn{x}{x}/2q} \Omega_x\Omega_{x/q}dx 
    \hspace{-43pt} & \\
  &= \left\langle\Omega_x,\Omega_x^{-1}\Omega_{x/p}\right\rangle_x
    \exp\left(
      S\left(\frac{p}{q}\right)
      -\frac{p}{q}-\frac{q}{p}-\frac{1}{pq}
      +3\sign(pq)
    \right)\frac{\theta}{48}
    & \parbox{92pt}{\raggedleft
      by Lemma~\ref{lem:SpecificIntegral} for $I(1,1/q;p/q)$
    } \\*
  &= \langle\Omega_x,\Omega_x^{-1}\Omega_{x/p}\rangle_x
    \exp\frac{-S(q/p)}{48}\theta
    & \text{by Definition~\ref{def:Dedekind}.}
\end{alignat*}
\end{proof}

\begin{corollary} \lbl{cor:SeparateLens}
The \LMO{} invariant does not separate lens spaces. (Though see
Section~\ref{subsec:SeifertIHS} for some more encouraging news about
the \LMO{} invariant).
\end{corollary}

\begin{proof} As noted in~\cite[pp.~247]{KirbyMelvin:DedekindSum},
the lens spaces $L_{25,4}$ and $L_{25,9}$ are not homeomorphic
but the Dedekind symbols $S(4/25)=S(9/25)=48/25$ are equal (see
Table~\ref{tab:STable}), and thus their \LMO{} invariants are equal.
\end{proof}

\begin{exercise} Recompute $\ZLMO(L_{p,q})$ directly from the Kontsevich
integral of Hopf chains (Proposition~\ref{prop:HopfChain}),
using the fact that surgery on a the Hopf chain
$H_{x_1,\ldots,x_\ell}(a_1,\ldots,a_\ell)$ gives $L_{q,-p}$.
\end{exercise}

\subsection{Seifert fibered spaces with a spherical base} Let
$M=S^3(b;p_1/q_1,\dots,p_n/q_n)$ be the Seifert fibered space with base
space $S^2$, with $n$ exceptional fibers with orbit invariants $(p_i,q_i)$
($p_i\geq 0$, $0<q_i<p_i$, $p_i$ and $q_i$ relatively prime), and with
bundle invariant $b\in{\mathbb Z}$. The following is easily established
(see~\cite{Montesinos:Tessellations, Scott:Geometries}):

\begin{itemize}
\item If $e_0:=b+\sum_i q_i/p_i\neq 0$, then $M$ is a rational homology
  sphere and $\left|H_1(M)\right|=|e_0|\prod_i p_i$.
\item $M$ is given by surgery on the following (rationally framed)
  ``key chain'' link in $S^3$:
  \[ L = L(b;p_1/q_1,\dots,p_n/q_n) := 
  \printname{KeyChain}
  \setlength{\unitlength}{0.5\standardunitlength}
  \begin{array}{c}  \hspace{-1.7mm}
    \raisebox{-8pt}{\begingroup\makeatletter\ifx\SetFigFont\undefined
\def\x#1#2#3#4#5#6#7\relax{\def\x{#1#2#3#4#5#6}}%
\expandafter\x\fmtname xxxxxx\relax \def\y{splain}%
\ifx\x\y   
\gdef\SetFigFont#1#2#3{%
  \ifnum #1<17\tiny\else \ifnum #1<20\small\else
  \ifnum #1<24\normalsize\else \ifnum #1<29\large\else
  \ifnum #1<34\Large\else \ifnum #1<41\LARGE\else
     \huge\fi\fi\fi\fi\fi\fi
  \csname #3\endcsname}%
\else
\gdef\SetFigFont#1#2#3{\begingroup
  \count@#1\relax \ifnum 25<\count@\count@25\fi
  \def\x{\endgroup\@setsize\SetFigFont{#2pt}}%
  \expandafter\x
    \csname \romannumeral\the\count@ pt\expandafter\endcsname
    \csname @\romannumeral\the\count@ pt\endcsname
  \csname #3\endcsname}%
\fi
\fi\endgroup
{\renewcommand{\dashlinestretch}{30}
\begin{picture}(8259,1956)(0,-10)
\put(1914.474,748.737){\arc{1297.477}{6.0810}{11.8412}}
\path(2592.618,762.881)(2550.000,879.000)(2532.959,756.486)
\put(3714.474,748.737){\arc{1297.477}{6.0810}{11.8412}}
\path(4392.618,762.881)(4350.000,879.000)(4332.959,756.486)
\put(6714.474,748.737){\arc{1297.477}{6.0810}{11.8412}}
\path(7392.618,762.881)(7350.000,879.000)(7332.959,756.486)
\put(750.000,1479.000){\arc{900.000}{1.5708}{4.7124}}
\put(7800.000,1479.000){\arc{900.000}{4.7124}{7.8540}}
\path(1500,1029)(3000,1029)
\path(3225,1029)(4725,1029)
\path(5550,1029)(6000,1029)
\path(6300,1029)(7800,1029)
\path(1200,1029)(750,1029)
\path(7800,1929)(4275,1929)
\path(4395.000,1959.000)(4275.000,1929.000)(4395.000,1899.000)
\path(4275,1929)(750,1929)
\put(450,1554){\makebox(0,0)[lb]{\smash{{\mathmode{\scriptstyle x}}}}}
\put(4950,579){\makebox(0,0)[lb]{\smash{{\mathmode{\cdots}}}}}
\put(0,1704){\makebox(0,0)[lb]{\smash{{\mathmode{\scriptstyle -b}}}}}
\put(1950,279){\makebox(0,0)[lb]{\smash{{\mathmode{\scriptstyle x_1}}}}}
\put(3750,279){\makebox(0,0)[lb]{\smash{{\mathmode{\scriptstyle x_2}}}}}
\put(6750,279){\makebox(0,0)[lb]{\smash{{\mathmode{\scriptstyle x_n}}}}}
\put(2250,54){\makebox(0,0)[lb]{\smash{{\mathmode{\scriptstyle p_1/q_1}}}}}
\put(4050,54){\makebox(0,0)[lb]{\smash{{\mathmode{\scriptstyle p_2/q_2}}}}}
\put(7050,54){\makebox(0,0)[lb]{\smash{{\mathmode{\scriptstyle p_n/q_n}}}}}
\end{picture}
} }
      \hspace{-1.9mm}
  \end{array}
 \]
\end{itemize}

We could use the rational surgery formula of Theorem~\ref{thm:main}
to compute the \LMO{} invariant of $M$, but it is somewhat
easier to backtrack a bit, and use the main points of the proof of
Theorem~\ref{thm:main}, applied in a slightly different manner. So we
shackle the components $x_1,\dots,x_n$ of $L$ as in Figure~\ref{fig:Ls},
and observe that the result $L^s$ is a $(-b)$-framed unknot marked $x$
with $n$ one-step-longer Hopf chains hanging from it. We then run the
``first step'' of Section~\ref{subsec:TwoStep} on these longer Hopf
chains. As there the result is $\Zw(\unknot^{-b}_x)$ corrected by $n$
insertions of factors corresponding to the $n$ one-step-longer Hopf
chains. So noting that $\langle 0,a_i^1,\dots,a_i^\ell\rangle = -1/\langle
a_i^1,\dots,a_i^\ell\rangle$ and using Equation~\eqref{eq:WheeledFormula}
we find that the result of the (revised) first step of the computation is

\begin{xalignat*}{1}
  \hat{\Omega}^{-1}\sigma Z'_1
  &= \Zw(\unknot_x^{-b})\prod_i \exp\left[
      \left(S(-q_i/p_i)+q_i/p_i\right)\frac{\theta}{48}
    \right]
    \Omega_x^{-1}\Omega_{x/p_i}
    \exp\left(-\frac{q_i}{2p_i}\strutn{x}{x}\right) \\
  &= \OmOm^{-1}
    \exp\left(\frac{\theta}{48}\left(
      e_0-\sum_i S\left(\frac{q_i}{p_i}\right)
    \right)\right)
    \Omega_x^{1-n}\left(\prod_i\Omega_{x/p_i}\right)
    \exp-\frac{e_0\strutn{x}{x}}{2}
\end{xalignat*}

\noindent
This result serves as the input into the second step of the computation in
the usual manner:

\begin{xalignat}{2}
  \hat{\Omega}^{-1}\sigma \Zc'_1
  &= \OmOm^{-2}    
    \exp\left(\frac{\theta}{48}\left( 
      e_0-\sum_i S\left(\frac{q_i}{p_i}\right)                     
    \right)\right)
    \Omega_x^{2-n}\left(\prod_i\Omega_{x/p_i}\right)
    \exp-\frac{e_0\strutn{x}{x}}{2} & \nonumber \\
  \Aa'_2
  &= \int \hat{\Omega}^{-1}\sigma \Zc'_1 \,dx & \nonumber \\
  &= \OmOm^{-2}    
    \exp\left(\frac{\theta}{48}\left( 
      e_0-\sum_i S\left(\frac{q_i}{p_i}\right)                     
    \right)\right)
    \left\langle
      \exp\frac{\strutu{x}{x}}{2e_0},
      \Omega_x^{2-n}\prod_i\Omega_{x/p_i}
    \right\rangle & \nonumber \\
  \ZLMO(M)
  &= \OmOm^{-1}
    \exp\left(\frac{\theta}{48}\left(
      \!e_0-3\sign(e_0)-\sum_i S\left(\frac{q_i}{p_i}\right)
    \right)\right)
    \left\langle
      \exp\frac{\strutu{x}{x}}{2e_0},
      \Omega_x^{2-n}\prod_i\Omega_{x/p_i}
    \right\rangle \lbl{eq:LMOSeifert}
\end{xalignat}

\begin{remark} \lbl{rem:CassonWalker}
According to~\cite[Proposition~5.3]{LeMurakamiOhtsuki:Universal}
(see also~\cite{LeMurakami2Ohtsuki:3ManifoldAnnouncement}), for any
rational homology sphere $M$ the coefficient of $\theta$ in $\ZLMO(M)$
is $\lambda_w(M)/4$, where $\lambda_w(M)$ denotes the Casson-Walker
invariant of $M$ (see~\cite{Walker:ExtensionOfCasson}). Thus in
our case, extracting the coefficient of $\theta$ in $\ZLMO(M)$ from
Equation~\eqref{eq:LMOSeifert}, we find that
\begin{equation} \lbl{eq:CassonWalker}
  \lambda_w(S^3(b;p_1/q_1,\dots,p_n/q_n)) = \frac{1}{12}\left(
     e_0 + \frac{2-n}{e_0} - 3\sign(e_0)
     + \sum_i\frac{1}{e_0p_i^2}-S\left(\frac{q_i}{p_i}\right)
  \right).
\end{equation}
This agrees with Lescop's evaluation of the same quantity
at~\cite[Proposition~6.1.1]{Lescop:CassonWalker} (though notice the
different normalization,~\cite[Section~1.5,~T5.0]{Lescop:CassonWalker}).
\end{remark}

\begin{remark} We do not know how to simplify the pairing in
Equation~\eqref{eq:LMOSeifert} any further, except in the case when
$n\leq 2$ (using Lemma~\ref{lem:SpecificIntegral}). But in this case $M$
is a lens space, so we learn nothing new. At any rate,
Remark~\ref{rem:CassonWalker} allows us to rewrite
Equation~\eqref{eq:LMOSeifert} in an alternate way:
\begin{equation} \lbl{eq:LMOSeifertAlt}
  \ZLMO(M) = \OmOm^{-1}\exp\frac{\lambda_w(M)\theta}{4}
    \exp\left(
      \frac{\theta}{48e_0}\left(n\!-\!2\!-\!\sum_i\frac{1}{p_i^2}\right)
    \!\right)
    \left\langle
      \exp\frac{\strutu{x}{x}}{2e_0},
      \Omega_x^{2-n}\prod_i\Omega_{x/p_i}
    \right\rangle
\end{equation}
\end{remark}

\begin{remark} \lbl{rem:GrossFailure}
Equation~\eqref{eq:LMOSeifertAlt} shows that the \LMO{} invariant
of Seifert fibered spaces with a spherical base factors through at
most $n+2$ integer parameters: the $p_i$'s, $e_0$ and the Casson-Walker
invariant of $M$.  As $S^3(b;p_1/q_1,\dots,p_n/q_n)$ depends on $2n$
integers (albeit with some range restrictions), it is clear that the
\LMO{} invariant grossly fails to separate rational homology Seifert
fibered spaces. Some more encouraging news about the \LMO{} invariant
are in Section~\ref{subsec:SeifertIHS}.
\end{remark}

\subsection{Some $sl(2)$ computations}

Recall that diagram spaces such as $\calA(\uparrow)$, $\calA(\ast)$
and $\calA(\emptyset)$ are in some sense ``universal''
versions of certain tensor spaces associated with Lie algebras
(see~\cite{Bar-Natan:Vassiliev}), and that the Kontsevich integral
and the \LMO{} invariant are in this sense universal versions
of the Reshetikhin-Turaev invariants of links and of manifolds
(see~\cite{Kassel:QuantumGroups, LeMurakami:Universal, Le:LMO,
Aarhus:I}). Thus we wish to explicitly compare our results with some
known answers for the Reshetikhin-Turaev invariants in the case of the
Lie algebra $sl(2)$.

In the particular case of $sl(2)$, we may divide all diagram spaces
$\calA(\ldots)$ by the additional ``$A_1$ relations'' (with $\hbar$
a formal parameter),
\begin{align}
  \bigcirc=3
  & \qquad\text{and}\qquad
  \HGraph=2\hbar\left(\hsmoothing-\crossing\right). \nonumber \\
\intertext{
  (The first relations says that $sl(2)$ is three dimensional,
  and the second becomes the identity $A\times(B\times C)=(A\cdot
  B)C-(A\cdot C)B$ in the vector calculus guise of $sl(2)$. See
  e.g.~\cite{ChmutovVarchenko:sl2}.) Calling the corresponding quotient
  spaces $\calA^1(\ldots)$, we easily find by repeated applications of
  the $A_1$ relations that
}
  \calA^1(\emptyset)=\bbQ[[\hbar]]
  & \qquad\text{and}\qquad
  \calA^1(\ast)=\bbQ[[\hbar,s]], \nonumber \\
\intertext{
  (here $s=(\frown)$ is a single strut).  One easily computes that in
  $\calA^1(\emptyset)$, $\theta=12\hbar$. Furthermore, simple induction
  (see e.g.~\cite{ThurstonD:Thesis}) shows that in $\calA^1(\ldots)$,
}
  \omega_{2m} = 2(2\hbar s)^m
  & \qquad\text{and}\qquad
  \langle s^m,s^m\rangle = (2m+1)! \lbl{eq:InA1}
\end{align}

We can now find in $\calA^1(\ldots)$:
\begin{xalignat*}{2}
   \Omega
  &= \exp_\udot \sum_{m=1}^\infty b_{2m}\omega_{2m}
    & \text{by Equation~\eqref{eq:OmegaDef}} \\
  &= \exp 2\sum_{m=1}^\infty b_{2m}(2\hbar s)^m
    & \text{by Equation~\eqref{eq:InA1}} \\
  &= \frac{\sinh\sqrt{\hbar s/2}}{\sqrt{\hbar s/2}}
    & \text{by Equation~\eqref{eq:bdef}} \\*
  &= \sum_{m=0}^\infty \frac{\hbar^m}{2^m(2m+1)!}s^m.
\end{xalignat*}
Hence by Equation~\eqref{eq:InA1},
\[
  \OmOm =\sum_{m=0}^\infty\frac{\hbar^{2m}}{2^{2m}(2m+1)!}
    = \frac{\sinh\hbar/2}{\hbar/2}.
\]
Assembling everything, we find from Equation~\eqref{eq:LMOSeifertAlt}
that in $\calA^1(\emptyset)$ the \LMO{} invariant of
$M=S^3(b;p_1/q_1,\dots,p_n/q_n)$ is given as a product of an easy factor
\[ \frac{\hbar/2}{\sinh\hbar/2} \exp3\hbar\lambda_w(M), \]
and a difficult `rest':
\begin{equation} \lbl{eq:RestAsBracket}
  Z^{\text{rest}}(M) := \left\langle
    \exp\frac{s}{2e_0},
    \frac
      {2\prod_i p_i\sinh\sqrt{\hbar s/2p_i^2}}
      {\hbar s\left(\sinh\sqrt{\hbar s/2}\right)^{n-2}}
  \right\rangle
  \exp\frac{\hbar}{4e_0}\left(n-2-\sum_i\frac{1}{p_i^2}\right).
\end{equation}

Next, we notice two equalities; the first is an immediate corollary of 
Equation~\eqref{eq:InA1}, and the second is an exercise in Gaussian
integration (for simplicity we restrict to $e_0>0$; otherwise we need to use
a different contour for the integration):
\[ 
  \left\langle\exp\frac{s}{2e_0},\left(\frac{\hbar s}{2}\right)^m\right\rangle
  = \frac{(2m+1)!\hbar^m}{(4e_0)^mm!}
  = \sqrt{\frac{\hbar^3e_0^3}{16\pi}}
    \int_\bbR d\beta\, e^{-\hbar e_0\beta^2/4}\beta^2(\hbar\beta/2)^{2m}.
\]
This known for every monomial in $\hbar s/2$, we can now replace the
bracket in Equation~\eqref{eq:RestAsBracket} by a Gaussian integral,
substituting $(\hbar\beta/2)^2$ for every occurrence of $\hbar s/2$
in~\eqref{eq:RestAsBracket}. Thus,
\begin{equation} \lbl{eq:RestAsIntegral}
   Z^{\text{rest}}(M) = \sqrt{\frac{e_0^3}{\pi\hbar}}
  \exp\frac{\hbar}{4e_0}\left(n-2-\sum_i\frac{1}{p_i^2}\right)
  \int_\bbR d\beta\,e^{-\hbar e_0\beta^2/4}\frac
    {\prod_i p_i\sinh\hbar\beta/2p_i}
    {\left(\sinh\hbar\beta/2\right)^{n-2}}.
\end{equation}
Disregarding some thorny normalization issues, this result is in agreement
with~\cite[Section~4.5]{LawrenceRozansky:SeifertManifolds}.

\subsection{The \LMO{} invariant of integral homology Seifert fibered
spaces} \lbl{subsec:SeifertIHS}

Corollary~\ref{cor:SeparateLens} and Remark~\ref{rem:GrossFailure}
seem like bad news for finite type invariants of 3-manifolds. But
things aren't as bad as they seem. The \LMO{} invariant is a universal
finite type invariant merely over the rationals, and it may well be that
rational homology spheres, whose homology groups may contain torsion,
are separated by finite-group-valued finite type invariants.  At any
rate, it is nice that we can complement our non-separation results for
the \LMO{} invariant of rational homology spheres with a separation
result for (some) integral homology spheres:

\begin{theorem} The \LMO{} invariant separates integral homology Seifert
fibered spaces.
\end{theorem}

\begin{proof} To be an integral homology sphere, a Seifert fibered
space must be of the form $M=S^3(b;p_1/q_1,\dots,p_n/q_n)$ discussed
above. Furthermore, it must have $\left|H_1(M)\right|=|e_0|\prod_i
p_i=|b+\sum_i q_i/p_i|\prod_i p_i=1$. For this to happen, $e_0$ must
equal $\pm 1\left/\prod_i p_i\right.$, the $p_i$'s must be pairwise
relatively prime, and then the $q_i$'s and $b$ are uniquely determined
up to an overall sign by the $p_i$'s, using the Chinese remainder
theorem. The overall sign of $b$ and the $q_i$'s (and therefore also
of $e_0$) can be read from the coefficient of $\theta$ in $\ZLMO(M)$
(see Equation~\eqref{eq:CassonWalker}). It remains to check to what
extent does $Z^{\text{rest}}(M)$ determine the $p_i$'s. Thus we regard
$Z^{\text{rest}}(M)$ as `known', and try to read out $n$ and
$p_1,\dots,p_n$.

For the sake of simplicity we assume that $e_0>0$;
that is, that $e_0=+1\left/\prod_i p_i\right.$. This done,
Equation~\eqref{eq:RestAsIntegral} shows that $Z^{\text{rest}}(M)$
makes sense as an honest analytic function of $\hbar>0$, and
not merely as a formal power series. Assume for the moment that
$n>2$. For large values of $\hbar$ it is easy to bound the integral
in Equation~\eqref{eq:RestAsIntegral} above and below by rational
functions of $\hbar$, and thus the growth rate of $Z^{\text{rest}}(M)$
as a function of $\hbar$ is determined by the exponential prefactor. Thus
by observing the growth rate of $Z^{\text{rest}}(M)$ we can determine
$(n-2-\sum_i 1/p_i^2)/4e_0$, factor that term out, and therefore regard
the integral in Equation~\eqref{eq:RestAsIntegral},
\[
  \sqrt{\frac{e_0^3}{\hbar}}
  \int_\bbR d\beta\,e^{-\hbar e_0\beta^2/4}\frac
    {\prod_i p_i\sinh\hbar\beta/2p_i}
    {\left(\sinh\hbar\beta/2\right)^{n-2}}
  = \frac{1}{\hbar^{3/2}}
  \int_\bbR d\beta\,e^{-\beta^2/4\hbar}\frac
    {\prod_i\sinh(\beta/2p_i\sqrt{e_0})}
    {\left(\sinh(\beta/2\sqrt{e_0})\right)^{n-2}}
\]
as known. As $\hbar$ varies, this is the integral of a fixed function
against all possible Gaussians. Knowing all of these integrals we know
what the function is. Thus we know all the quantities $p_i\sqrt{e_0}$
and the value of $\sqrt{e_0}$ (recall that the $p_i$'s are all distinct
and greater than $1$, so no accidental cancellations can occur). This
finishes the case of $n>2$.

If $n\leq 2$, then $M$ is a sphere and therefore $\ZLMO(M)=1$. This
leads to $Z^{\text{rest}}(M)=\frac{\sinh\hbar/2}{\hbar/2}$, which
grows like $e^{\hbar/2}$. This growth rate is much slower than the
minimal possible value of $(n-2-\sum_i 1/p_i^2)/4e_0$ for $n>2$, which
is attained when $n=3$ and $(p_i)=(2,3,5)$. Thus the case of $n\leq 2$
is easily separated from the case of $n>2$.
\end{proof}

\begin{exercise} Verify directly from Equation~\eqref{eq:RestAsIntegral}
that in the $sl(2)$ case $Z^{\text{rest}}(M)=\frac{\sinh\hbar/2}{\hbar/2}$
for $n\leq 2$.
\end{exercise}

\ifx\undefined\bysame
        \newcommand{\bysame}{\leavevmode\hbox to3em{\hrulefill}\,}
\fi

\parpic(1.8in,1.5in)[l]{\includegraphics{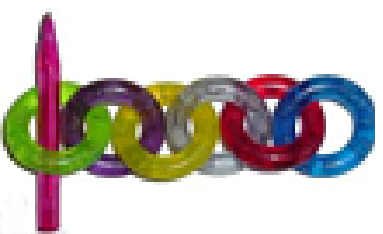}}
 

\begin{thebibliography}{LMMO}

\bibitem[B-N]{Bar-Natan:Vassiliev} D.~Bar-Natan,
  {\em On the Vassiliev knot invariants,}
  Topology {\bf 34} 423--472 (1995).

\bibitem[BGRT]{Bar-NatanGaroufalidisRozanskyThurston:WheelsWheeling}
  \bysame, S.~Garoufalidis, L.~Rozansky and D.~P.~Thurston,
  {\em Wheels, wheeling, and the Kontsevich integral of the unknot,}
  to appear in the Israel Jour.{} of Math, arXiv:q-alg/9703025.

\bibitem[BLT]{Bar-NatanLeThurston:InPreparation} \bysame, T.~Q.~T.~Le,
  and D.~P.~Thurston, in preparation.

\bibitem[CV]{ChmutovVarchenko:sl2} S.~V.~Chmutov and A.~N.~Varchenko,
  {\em Remarks on the Vassilliev knot invariants coming from $sl_2$,}
  Topology {\bf 36-1} (1997).

\bibitem[HV]{HinichVaintrob:CyclicOperads} V.~Hinich and A.~Vaintrob,
  {\em Cyclic operads and algebra of chord diagrams,}
  University of Haifa and IHES preprint, May 2000, arXiv:math.QA/0005197.

\bibitem[Kas]{Kassel:QuantumGroups} C.~Kassel,
  {\em Quantum groups,}
  Springer-Verlag GTM {\bf 155}, Heidelberg 1994.

\bibitem[KM]{KirbyMelvin:DedekindSum} R.~Kirby and P.~Melvin,
  {\em Dedekind sums, $\mu$-invariants and the signature cocycle,}
  Math.{} Ann.{} {\bf 299} (1994) 231--267.

\bibitem[Ko]{Kontsevich:DeformationQuantization} M.~Kontsevich,
  {\em Deformation quantization of Poisson manifolds,}
  I.H.E.S.\ preprint, September 1997, arXiv:q-alg/9709040.

\bibitem[LR]{LawrenceRozansky:SeifertManifolds} R.~Lawrence and L.~Rozansky,
  {\em Witten-Reshetikhin-Turaev invariants of Seifert manifolds,}
  Comm.{} Math.{} Phys.{} {\bf 205} (1999) 287--314.

\bibitem[Le]{Le:LMO} T.~Q.~T.~Le,
  {\em The \LMO{} invariant,}
  SUNY at Buffalo preprint, June 1999.

\bibitem[LM1]{LeMurakami:Universal} \bysame\ and J.~Murakami,
  {\em The universal Vassiliev-Kontsevich invariant for framed oriented links,}
  Compositio Math.\ {\bf 102} (1996), 41--64, arXiv:hep-th/9401016.

\bibitem[LM2]{LeMurakami:Parallel} \bysame\ and \bysame,
  {\em Parallel version of the universal Vassiliev-Kontsevich invariant,}
  J.\ Pure and Appl.\ Algebra {\bf 121} (1997) 271--291.

\bibitem[LMMO]{LeMurakami2Ohtsuki:3ManifoldAnnouncement} \bysame,
  \bysame, H.~Murakami and T.~Ohtsuki,
  {\em A three-manifold invariant derived from the universal
    Vassiliev-Kontsevich invariant},
  Proc.\ Japan Acad.\ {\bf 71} Ser.\ A (1995) 125--127.

\bibitem[LMO]{LeMurakamiOhtsuki:Universal} \bysame, \bysame\ and
  T.~Ohtsuki,
  {\em On a universal perturbative invariant of 3-manifolds,}
  Topology {\bf 37-3} (1998) 539--574, arXiv:q-alg/9512002.

\bibitem[Les]{Lescop:CassonWalker} C.~Lescop,
  {\em Global surgery formula for the Casson-Walker invariant,}
  Annals of Mathematics Studies {\bf 140}, Princeton University Press,
  Princeton 1996.

\bibitem[Me]{Mermin:Moon} N.~D.~Mermin,
  {\em Is the moon there when nobody looks? Reality and the quantum
    theory,}
  Phys.\ Today {\bf 39-4} (1985) 38--47.

\bibitem[Moc]{Mochizuki:KirillovDuflo} T.~Mochizuki,
  {\em On the Kirillov-Duflo type morphism,}
  Osaka City University preprint, June 1999.

\bibitem[Mon]{Montesinos:Tessellations} J.~M.~Montesinos,
  {\em Classical tessellations and three-manifolds,}
  Springer-Verlag, 1985.

\bibitem[Ro]{Rolfsen:KnotsLinks} D.~Rolfsen,
  {\em Knots and Links},
  Publish or Perish, Mathematics Lecture Series {\bf 7}, Wilmington
  1976.

\bibitem[Sc]{Scott:Geometries} P.~Scott,
  {\em The geometries of 3-manifolds,}
  Bull.{} London Math.{} Soc.{} {\bf 15} (1983) 401--487.

\bibitem[Th]{ThurstonD:Thesis} D.~P.~Thurston,
  {\em Wheeling: A Diagrammatic Analogue of the Duflo Isomorphism,}
  University of California at Berkeley Ph.D.{} thesis, May 2000,
  arXiv:math.QA/0006083.

\bibitem[Wa]{Walker:ExtensionOfCasson} K.~Walker,
  {\em An extension of Casson's invariant,}
  Annals of Mathematics Studies {\bf 126}, Princeton University Press,
  Princeton 1992.

\bibitem[\AA-I]{Aarhus:I}
  D.~Bar-Natan, S.~Garoufalidis, L.~Rozansky and D.~P.~Thurston,
  {\em The \AA{}rhus integral of rational homology 3-spheres I: A highly
    non trivial flat connection on $S^3$,}
  Selecta Math., to appear, arXiv:q-alg/9706004.

\bibitem[\AA-II]{Aarhus:II}
  \bysame, \bysame, \bysame\ and \bysame,
  {\em The \AA{}rhus integral of rational homology 3-spheres II:
        Invariance and Universality,}
  Selecta Math., to appear, arXiv:math.QA/9801049.

\bibitem[\AA-III]{Aarhus:III}
  \bysame, \bysame, \bysame\ and \bysame,
  {\em The \AA{}rhus integral of rational homology 3-spheres III:
    The relation with the Le-Murakami-Ohtsuki invariant,}
  Hebrew University, Harvard University, University of Illinois at
  Chicago and University of California at Berkeley preprint, August 1998,
  arXiv:math.QA/9808013.

\end{thebibliography}
\end{document}